\documentclass[11pt,a4paper]{amsart}
\usepackage{amsmath,amssymb,amsthm,mathtools,esint,bm}
\usepackage[hidelinks]{hyperref}
\usepackage{enumitem,microtype,xcolor}

\newcommand{\Scone}{\mathcal S}
\newcommand{\Hcone}{\mathcal H}
\newcommand{\del}{\delta_D}
\newcommand{\esssup}{\operatorname*{ess\,sup}}
\newcommand{\essinf}{\operatorname*{ess\,inf}}
\newcommand{\osc}{\operatorname{osc}}
\theoremstyle{plain}
\newtheorem{theorem}{Theorem}[section]
\newtheorem{lemma}[theorem]{Lemma}
\newtheorem{proposition}[theorem]{Proposition}
\newtheorem{corollary}[theorem]{Corollary}
\theoremstyle{definition}
\newtheorem{definition}[theorem]{Definition}
\newtheorem{remark}[theorem]{Remark}
\theoremstyle{remark}
\newtheorem{example}[theorem]{Example}
\numberwithin{equation}{section}
\title[Boundary Harnack principle for abstract supersolutions]{{Boundary Harnack principle for abstract supersolutions on uniform domains}}
\author{Takanobu Hara}
\date{\today}
\begin{document}

\begin{abstract}
We develop an axiomatic approach to the boundary Harnack principle on diameter-uniform domains in proper doubling metric measure spaces. Under critical-density and upper-contraction assumptions, we establish a scale-invariant one-sided boundary comparison between nonnegative supersolutions and positive solutions as well as H\"{o}lder continuity of quotients of positive solutions. Our framework does not require linearity or rectifiable curves. Applications include classical elliptic equations with Dirichlet and Neumann boundary conditions, spectral fractional Laplacians, harmonic functions on fractals, and Pucci operators.
\end{abstract}
\maketitle

\section{Introduction}\label{sec:introduction}

The boundary Harnack principle (BHP) compares positive harmonic functions that vanish on a common portion of the boundary. For Lipschitz domains, fundamental results were obtained by Ancona \cite{An78}, Dahlberg \cite{D77}, and Wu \cite{W78}; Caffarelli, Fabes, Mortola, and Salsa \cite{CFMS81} established boundary comparison for divergence-form elliptic operators. Jerison and Kenig \cite{JK82} proved the BHP on non-tangentially accessible domains and obtained H\"older regularity of quotients of positive harmonic functions \cite[Theorem~7.9]{JK82}. Their domain condition includes exterior accessibility as well as interior geometry. Bass and Burdzy \cite{BB91} proved a global BHP on twisted H\"older domains by probabilistic methods without imposing exterior regularity.

For uniform domains, Aikawa \cite{A01} proved a scale-invariant BHP for harmonic functions, together with H\"older regularity of quotients and applications to the Martin boundary. Gyrya and Saloff-Coste \cite{GS11} extended Aikawa's approach to symmetric strongly local Dirichlet spaces and obtained the BHP also for inner uniform domains. Lierl and Saloff-Coste \cite{LSC14} proved a scale-invariant BHP on inner uniform domains for a broader class of local Dirichlet forms, including non-symmetric ones. Chen \cite{Chen25} established boundary comparison on uniform domains in the setting of strongly local symmetric Dirichlet forms.

{\color{black}Our approach to the BHP on uniform domains relies on the following two ingredients.}
\begin{enumerate}[label=(\arabic*)]
\item \emph{Geometry of uniform domains.} Gehring and Osgood \cite{GO79} related uniform domains to quasihyperbolic distance estimates. Aikawa and Shanmugalingam \cite{AS05} studied reference-point constructions and boundary-layer estimates in the geometry of uniform and John domains. The metric formulation of Eriksson-Bique and Gong \cite{EBG21} makes it possible to work with subcurve diameters without requiring rectifiable curves.
\item \emph{An axiomatic framework for Harnack inequalities.} Aimar, Forzani, and Toledano \cite{AFT} developed an abstract approach to Harnack inequalities based on properties of classes of functions. Safonov \cite{Safonov} used growth lemmas to prove boundary comparison directly, without relying on Green functions or harmonic measure, and later extended this approach to metric spaces \cite{Saf20}. These works provide a basis for studying boundary comparison through properties of solutions and domain geometry.
{For further historical background, see \cite{Kas07,Safonov,Saf20}.}
\end{enumerate}

Our main results are a scale-invariant one-sided boundary comparison for nonnegative supersolutions and positive solutions (Theorem~\ref{thm:bhp}),
and a two-point H\"{o}lder estimate for the quotient of two positive solutions (Theorem~\ref{thm:quotient}).
The framework is based on a comparison-difference property (Comp), a critical-density estimate (CD), and an upper-contraction property (UC).
The latter is also required near the boundary for suitable negative parts of comparison differences.
Boundary (UC) is a condition on comparison functions rather than a prescribed boundary condition; its Dirichlet and Neumann realizations are treated in Section~\ref{sec:verifyA2}.

Our proof follows the PDE-based approaches of Safonov \cite{Safonov}, Kim and Safonov \cite{KS11, KS14}, De Silva and Savin \cite{DS20}, and Maiale, Tortone, and Velichkov \cite{MTV24}. Building on these PDE-based approaches, we incorporate quantitative estimates developed in potential theory, specifically localized ball-chain bounds and boundary-layer measure decay derived from interior porosity. This allows us to exploit the geometry of the domain more fully and extend the boundary Harnack principle to a broader setting.

The present axioms are tailored to boundary comparison and do not impose linearity, rectifiable curves,
or certain additional structural properties, such as locality and closure under minima.
As shown in Section~\ref{sec:examples}, the framework applies to divergence- and nondivergence-form elliptic equations, spectral fractional Laplacians,
harmonic functions on the Sierpi\'{n}ski gasket, and fully nonlinear Pucci operators.
{{\color{black}The nonlinear equations covered by the present framework admit estimates for comparison differences.} For nonlinear variational equations such as the $p$-Laplace equation, a different approach based on energy estimates and the Bregman divergence is developed in \cite{HaraPH}. Although used in different ways, logarithmic bounds on the lengths of Harnack chains and boundary porosity are the key geometric ingredients in both approaches.}
On the other hand, Section~\ref{sec:limits} illustrates the limitations of these axioms, and their implications for Martin and Kuramochi boundary theory remain to be investigated.

Section~\ref{sec:ambient} states the ambient and domain assumptions, and Section~\ref{sec:geometry} proves the geometric estimates. Section~\ref{sec:interior} derives the interior consequences of the density axioms. Section~\ref{sec:boundary} proves boundary comparison and quotient regularity. Sections~\ref{sec:examples} and~\ref{sec:limits} treat realizations and limitations of the assumptions. {Section~\ref{sec:verifyA2} gives criteria for verifying (UC) in concrete equations.}

\section{Ambient setting and diameter uniformity}\label{sec:ambient}

Let $(X,d,\mu)$ be a proper metric measure space. There is a constant $C_\mu\ge1$ such that, for every $x\in X$ and $s>0$, the Borel measure $\mu$ satisfies
\begin{equation}\label{eq:doubling}
0<\mu(B(x,2s))\le C_\mu\mu(B(x,s))<\infty.
\end{equation}
Write $Q=\log_2C_\mu$.
Let $D\subsetneq X$ be nonempty and open, and set $\del(x)=d(x,X\setminus D)$.

\paragraph{Diameter-uniform domain.}
There is $A\ge1$ such that every $x,y\in D$ can be joined by a continuous curve $\gamma:[0,1]\to D$ satisfying, for every $t\in[0,1]$,
\begin{equation}\label{eq:diamuniform}
 \operatorname{diam}\gamma\le Ad(x,y),\qquad
 \del(\gamma(t))\ge A^{-1}\min\{\operatorname{diam}\gamma([0,t]),
 \operatorname{diam}\gamma([t,1])\}.
\end{equation}
We use the subcurve-diameter convention of \cite[Definition~4.12]{EBG21}, rather than the endpoint-distance convention of \cite[Definition~3.1]{Chen25}. The relation to length-uniformity, including equivalence in doubling quasiconvex ambient spaces, is discussed following \cite[Definition~4.12]{EBG21}.

\begin{example}[Snowflaked uniform domains]\label{ex:snowflake}
The standard snowflake metric $d_\theta(x,y)=|x-y|^\theta$, $0<\theta<1$, illustrates the distinction. The identity from Euclidean space to $(\mathbb R^n,d_\theta)$ is quasisymmetric with distortion $t\mapsto t^\theta$, so it preserves diameter-uniformity by \cite[Remark~4.16(4)]{EBG21}. For a Euclidean length-uniform domain with constant $A_E$, the same curves satisfy~\eqref{eq:diamuniform} with constant $A_E^\theta$. Yet no nonconstant curve is rectifiable for $d_\theta$ \cite[Section~1, p.~1]{KLD14}, so the domain is not length-uniform in that metric.

With Lebesgue measure, $(\mathbb R^n,d_\theta)$ is proper, connected, unbounded, and doubling.
\end{example}

{Fix}
\[
 r_* = \frac{\operatorname{diam}D}{16A},\qquad \kappa=\frac{1}{8A},
\]
with the convention $r_*=\infty$ when $\operatorname{diam}D=\infty$. All boundary scales below are smaller than $r_*$. Geometric constants depend only on $A$ and $C_\mu$ (and any prescribed ball dilation). In the analytic statements, structural constants may additionally depend on the common constant $M$ in (CD) and (UC) and the fixed analytic dilation factors. {Set $\Lambda=16(A+2)$. Lemma~\ref{lem:chains} uses balls with radii chosen after the required dilation factor has been fixed; their dilates then lie in $D\cap B(\xi,\Lambda s)$. The boundary scale factor $\Lambda_0$ is fixed in Section~\ref{sec:boundary}, after the interior Harnack dilation $C_{\mathrm{int}}$ has been obtained.}

Measure doubling gives metric doubling and the packing bounds used below, by comparing disjoint small balls with a containing ball.

\section{Geometric consequences of diameter uniformity}\label{sec:geometry}
The constructions of interior reference points and boundary layers adapt \cite[Remark~4.2 and proof of Lemma~4.4]{AS05} to continuous curves and packing in the present setting.

Every open subset $U$ of $D$ is locally path connected. Indeed, if $B(x,R)\subset U$ and $d(x,y)<R/(2A)$, a curve from~\eqref{eq:diamuniform} joins $x$ to $y$ inside $B(x,R)$. Thus the path components of $U$ are open, and every connected open subset of $D$ is path connected.

\begin{lemma}[Interior anchors]\label{lem:anchors}
{Let $\xi\in X\setminus D$ and $0<s<r_*$ satisfy $d(\xi,D)<s/16$. Then there exists $a_{\xi,s}\in D\cap B(\xi,s)$ such that $\del(a_{\xi,s})\ge\kappa s$.}
\end{lemma}
\begin{proof}
{Choose $x\in D\cap B(\xi,s/16)$ and $y\in D$ with $d(y,\xi)>2s$,} using $\operatorname{diam}D>4s$ (also valid when the diameter is infinite). On a uniform curve from $x$ to $y$, continuity provides $a$ with $d(a,\xi)=s/2$. The two subcurve diameters at $a$ are at least $d(x,a)>7s/16$ and $d(y,a)>3s/2$. {Hence,} Equation~\eqref{eq:diamuniform} gives $\del(a)>7s/(16A)\ge\kappa s$.
\end{proof}

\begin{lemma}[Chains with controlled center distances]\label{lem:packingchain}
{Suppose a continuous curve $\sigma\subset D$ has diameter at most $Kb$ and $\del\ge b$ on it. Fix $\lambda\ge1$, $0<\theta<1$, and $0<c<\min\{1/4,1/(2\lambda)\}$. Its endpoints can be compared through at most $C(K,c,\theta,C_\mu)$ interior balls of radius $cb$. Every consecutive pair of centers can be taken at distance smaller than $\theta cb$, and endpoint-to-center distances satisfy the same bound. The centers lie on $\sigma$, and the balls $B(z,\lambda cb)$ lie in $D$.}
\end{lemma}
\begin{proof}
Take a maximal $\varepsilon$-separated subset of the compact curve image, where $\varepsilon=\theta cb/8$. Metric doubling bounds its cardinality by $C(1+K/(\theta c))^Q$. The balls of radius $2\varepsilon$ cover the image. Their relative intersections with that image have a connected intersection graph. Choose a simple path in this graph between balls containing the endpoints, discarding unused balls. Consecutive centers have distance less than $4\varepsilon<\theta cb$. {Since $d(z,X\setminus D)\ge b$ and $\lambda cb<b/2$, each ball $B(z,\lambda cb)$ lies in $D$.}
\end{proof}

\begin{lemma}[Localized logarithmic and short chains]\label{lem:chains}
{Fix $\lambda\ge1$. There is $C\ge1$, depending only on $A$, $C_\mu$, and $\lambda$, with the following properties.}
\begin{enumerate}[label=(\roman*)]
\item {If $\xi\in X\setminus D$, $0<s<r_*$, $d(\xi,D)<s/16$, and $x\in D\cap B(\xi,4s)$, then there is a chain to $a=a_{\xi,s}$ with all $\lambda$-dilated balls in $D\cap B(\xi,\Lambda s)$, using at most}
\begin{equation}\label{eq:chainN}
 N(x,a)\le C\left(1+\log^+\frac{s}{\del(x)}\right)
\end{equation}
balls.
\item {If $x\in D$, $t/2\le\del(x)<t$, and $0<t<r_*$, then $x$ can be joined to a ball in $\{\del\ge t\}$ by at most $C$ balls, all their $\lambda$-dilates lying in $D\cap B(x,Ct)$.}
\end{enumerate}
\end{lemma}
\begin{proof}
(i) Choose a diameter-uniform curve from $x$ to $a$, whose diameter is at most $5As$. The functions $t\mapsto\operatorname{diam}\gamma([0,t])$ and $t\mapsto\operatorname{diam}\gamma([t,1])$ are continuous, respectively nondecreasing and nonincreasing. Split at a point where they agree, with common value $R\le5As$. On the half starting at $x$, the cigar condition gives $\del\ge A^{-1}\operatorname{diam}\gamma([0,t])$. Before the prefix diameter reaches $\del(x)/4$, the Lipschitz property of $\del$ instead gives $\del\ge3\del(x)/4$. Divide the rest by successive prefix diameters $2^j\del(x)/4$, stopping at $R$. Each resulting subcurve has diameter at most twice its lower scale and boundary distance at least that lower scale divided by $A$. Lemma~\ref{lem:packingchain} gives a bounded number of balls per scale. If the initial scale exceeds $R$, the entire half has controlled diameter relative to $\del(x)$ and needs only a bounded number of balls. The same construction from $a$ has a bounded number of scales because $\del(a)\ge\kappa s$. Thus~\eqref{eq:chainN} follows.

For each subcurve, adjoin balls centered at its two endpoints, with the same radius as the balls supplied by Lemma~\ref{lem:packingchain}. The endpoint-to-center bound in that lemma preserves the controlled center distances. The radii associated with two subcurves meeting at a junction are comparable. Within either half, successive dyadic scales differ by a factor of two, and the initial segment has a radius comparable to the first dyadic radius. At the splitting point, the terminal scales on the two halves are comparable to $R$ whenever both halves have a dyadic segment. If one half consists only of its initial segment, its endpoint boundary distance exceeds a fixed multiple of $R$; the $1$-Lipschitz property of $\del$ and $\operatorname{diam}\gamma\le2R$ then show that the endpoint boundary distances, and hence the two terminal radii, are comparable. Consequently, for adjacent radii $r_j,r_{j+1}$ there is $C=C(A)\ge1$ such that $C^{-1}\le r_j/r_{j+1}\le C$. At a junction, use the common endpoint as a center and, if necessary, insert a ball with radius $\min\{r_j,r_{j+1}\}$. The two chains can then be connected by a uniformly bounded number of comparisons in balls of the smaller radius. Choose the radius factor $c$ in Lemma~\ref{lem:packingchain} after fixing the desired ball dilation $\lambda$ (in Harnack applications, at least the interior Harnack dilation), and take $\theta<1/4$. Each center comparison is then performed in a ball whose required dilation is contained in $D$. Finally, the curve lies in $B(\xi,(4+5A)s)$ and each center $z$ on it satisfies $\del(z)\le d(z,\xi)\le(4+5A)s$. Taking $c$ sufficiently small also places every required dilated ball in $B(\xi,\Lambda s)$, where $\Lambda=16(A+2)$. Thus all such balls lie in $D\cap B(\xi,\Lambda s)$.

(ii) Choose $y\in D$ with $d(x,y)>8At$ and a uniform curve from $x$ to $y$. Such a point exists since $t<r_*$: for finite diameter, $\sup_{y\in D}d(x,y)\ge\operatorname{diam}D/2>8At$, and for infinite diameter the distances from $x$ are unbounded. At the first prefix diameter $2At$, the remaining subcurve has diameter greater than $2At$. Thus its boundary distance is at least $2t$. Along the preceding part, $\del\ge t/(4A)$: before prefix diameter $t/4$ this follows from $\del(x)\ge t/2$; subsequently it follows from the cigar condition, since the remaining diameter is greater than $2At$. This part has diameter $2At$. {Hence, Lemma~\ref{lem:packingchain} gives a uniformly bounded number of balls of radius a sufficiently small multiple of $t$.} Its terminal ball has $\del\ge t$. {The $\lambda$-dilated balls lie in $D\cap B(x,Ct)$ by the choice of $c$ in Lemma~\ref{lem:packingchain}.}
\end{proof}

\begin{lemma}[Boundary layers]\label{lem:layer}
{There are structural $\alpha_{\mathrm{lay}}>0$ and $C_L\ge1$ such that, for every $\xi\in X\setminus D$, $0<s<r_*$ with $d(\xi,D)<s/16$, and $0<t<s$,}
\begin{equation}\label{eq:layer}
 \frac{\mu(D\cap B(\xi,4s)\cap\{\del<t\})}{\mu(B(\xi,4s))}
 \le C_L(t/s)^{\alpha_{\mathrm{lay}}}.
\end{equation}
{Consequently, for every $0<\eta<\alpha_{\mathrm{lay}}$ there exists $C_\eta<\infty$, depending only on $A$, $C_\mu$, and $\eta$, such that}
\begin{equation}\label{eq:localnegativepower}
 \int_{D\cap B(\xi,4s)}(s/\del(x))^\eta\,d\mu(x)
 \le C_\eta\mu(B(\xi,4s)).
\end{equation}
\end{lemma}
\begin{proof}
Fix $a=a_{\xi,s}$. {Set} $t_0=\kappa s/(A+2)$, $t_j=2^{-j}t_0$ and
\[
 V_j=\{x\in D\cap B(\xi,4s+{8(A+2)}t_j):t_j/2\le\del(x)<t_j\},\qquad E_j=\bigcup_{i>j}V_i.
\]
For $x\in E_j$, take a uniform curve to $a$ and let $y$ be its first point with $\del(y)=3t_j/4$. This point exists by continuity. The remaining subcurve has diameter at least $\del(a)-3t_j/4>3At_j/4$. The cigar condition at $y$ therefore forces its prefix diameter to be at most $3At_j/4$. Hence $d(x,y)\le3At_j/4$. The ball $B(y,t_j/8)$ has boundary distances in $(t_j/2,t_j)$ and lies in $V_j$, since $x\in B(\xi,4s+{4(A+2)}t_j)$ and ${4(A+2)+3A/4+1/8<8(A+2)}$. It also lies in $B(x,(A+2)t_j)$. Apply the $5$-covering lemma to the latter balls. The selected balls are disjoint, as are their contained small balls. Doubling compares the measures of the fivefold selected balls with those of the small balls in $V_j$. Thus
\[
 \sum_{i>j}\mu(V_i)=\mu(E_j)\le C\mu(V_j).
\]
It follows that $S_j=\sum_{i\ge j}\mu(V_i)$ satisfies $S_{j+1}\le qS_j$ with a structural $q<1$. Also $S_0\le\mu(B(\xi,5s))$. The thin set with $\del<t_j$ in $B(\xi,4s)$ lies in $\bigcup_{i\ge j}V_i$. Choose $\alpha_{\mathrm{lay}}>0$ with $2^{\alpha_{\mathrm{lay}}} q<1$. {Thus, the geometric decay of $S_j$ and comparison with consecutive dyadic scales prove~\eqref{eq:layer} for $t<t_0$. Increasing $C_L$ covers $t\ge t_0$. Summing over dyadic boundary layers then yields~\eqref{eq:localnegativepower}.}
\end{proof}

\begin{remark}[Relation to the length formulation]\label{rem:lengthrelation}
{Suppose in addition that $D$ is length-uniform with constant $A_\ell$: any two points $x,y\in D$ can be joined by a rectifiable curve of length at most $A_\ell d(x,y)$ satisfying the length cigar condition with the same constant. Let $\xi\in X\setminus D$, $0<s<r_*$, and $d(\xi,D)<s/16$. For $x\in D\cap B(\xi,4s)$ set $a=a_{\xi,s}$. Let $k_{\xi,s}(x,a)$ be the infimum of $\int_\gamma\del(z)^{-1}\,ds(z)$ over rectifiable curves in $D\cap B(\xi,(4+5A_\ell)s)$ joining $x$ to $a$. The length-uniform curve from $x$ to $a$ has length at most $5A_\ell s$ and lies in this neighborhood. Integrating along its endpoint halves gives}
\[
 k_{\xi,s}(x,a)\le C(1+\log^+(s/\del(x))),
\]
{where $C$ depends only on $A$ and $A_\ell$. Combining this estimate with~\eqref{eq:localnegativepower} yields the corresponding exponential quasihyperbolic moment for every positive exponent smaller than $\alpha_{\mathrm{lay}}/C$.} For Euclidean H\"older domains, exponential quasihyperbolic integrability is established in \cite[Theorem~A]{SS91}.
\end{remark}

\section{{Supersolution families} and interior estimates}\label{sec:interior}\label{sec:analytic}

 The abstract Harnack framework of Aimar, Forzani, and Toledano \cite{AFT} guides the density axioms below. To obtain their consequences, we adapt the density and covering arguments in \cite[Sections~4--7]{AFT} and the proofs in \cite[Theorems~4.7 and~5.1]{DFGL08}. For the interior Harnack inequality, we follow the weighted-maximum proof in \cite[Theorem~3.1]{Safonov}. {The arguments use essential lower bounds for supersolutions and continuous representatives for positive solutions; (Comp) does not require negative constants to belong to $\mathcal P$.}

The interior estimates use the ambient setting and (Comp), (CD), and (UC) on interior balls, without a boundary assumption. Diameter uniformity enters only through local path connectedness in the strong minimum principle.

\begin{definition}[(Comp) Comparison differences]\label{def:Comp}
For every open set $U\subset D$, let $\Scone(U)$, $\Hcone(U)$,
and $\mathcal P(U)$ be classes of real-valued measurable functions on
$U$, compatible with restriction. Assume that
$0\in\Hcone(U)\subset\Scone(U)$ and that positive members of
$\Hcone(U)$ are continuous.

For every $u\in\Scone(U)$, $v\in\Hcone(U)$, and $t\ge0$, require
\[
 tu-v\in\mathcal P(U).
\]
The class $\mathcal P(U)$ is closed under multiplication by positive
scalars and addition of nonnegative constants. Define
\[
 \mathcal C(U)=\{w\in\mathcal P(U):-w\in\mathcal P(U)\}.
\]
\end{definition}
For every open $U\subset D$, (Comp) implies
$\Scone(U)\subset\mathcal P(U)$ and
$\Hcone(U)\subset\mathcal C(U)$. If $u,v\in\Hcone(U)$ and
$t\ge0$, then $tu-v\in\mathcal C(U)$: for $t>0$,
$v-tu=t(t^{-1}v-u)\in\mathcal P(U)$, and the case $t=0$
follows from $\Hcone(U)\subset\mathcal C(U)$.
Members of $\Hcone(U)$, $\Scone(U)$, and $\mathcal P(U)$
are called solutions, supersolutions, and comparison supersolutions,
respectively.

\begin{definition}[{(CD) Critical-density property}]\label{def:CD} Put $\vartheta=(2C_\mu^2)^{-1}$ and fix $A_{\mathrm{CD}}\ge4$. There is $M\ge2$ such that {for every open $U\subset D$, every nonnegative $u\in\mathcal P(U)$, every $x\in X$ and $s>0$ satisfying $B(x,A_{\mathrm{CD}}s)\subset U$, and every $t>0$,}
\begin{equation}\label{eq:CD}
 \frac{\mu(B(x,s)\cap\{{u}\ge t\})}{\mu(B(x,s))}\ge\vartheta
 \quad\Longrightarrow\quad \essinf_{B(x,s/2)}{u}\ge t/M.
\end{equation}
\end{definition}
The density threshold is determined by ambient doubling, not by a dimension. Applying (CD) on concentric balls yields the double-ball property in \cite[Section~2]{AFT} and \cite[Definition~4.2]{DFGL08}; its quantitative propagation form is proved in Lemma~\ref{lem:propagation}. This property of comparison supersolutions differs from measure doubling~\eqref{eq:doubling}.

\begin{definition}[{(UC)}]\label{def:UC}
Let $M\ge2$ be the constant in (CD), enlarged if necessary. Let
$V\subset X$ be open and let ${f}\ge0$ be defined on $D\cap V$.
We say that (UC) holds for ${f}$ on $V$ if ${f}$ is essentially
bounded on $D\cap B(x,s)$ for every ball $B(x,2s)\subset V$ and,
whenever $B(x,4s)\subset V$ and $t>0$,
\begin{equation}\label{eq:UD}
 \frac{\mu(D\cap B(x,2s)\cap\{{f}>t\})}{\mu(B(x,2s))}\le M^{-1}
 \quad\Longrightarrow\quad
 \esssup_{D\cap B(x,s)}{f}\le Mt+{\frac{1}{2}}\esssup_{D\cap B(x,2s)}{f}.
\end{equation}
The essential supremum of a nonnegative function over the empty set
is zero. When (UC) is specified only for balls centered in $D$, both the
boundedness requirement and \eqref{eq:UD} are restricted to such centers;
otherwise all centers in $X$ are permitted. For every open $U\subset D$,
we assume (UC) on $U$ for ${f}=({u}-c)_+$ whenever ${u}$ is a nonnegative
continuous member of $\mathcal C(U)$ and $c\ge0$. Boundary applications require (UC) separately
for the functions and centers specified in their statements. The constant
$M$ is uniform, and no value of ${f}$ on $X\setminus D$ is involved.
The boundary arguments apply (UC) at centers chosen in $D$; the outer
ambient balls may nevertheless meet $X\setminus D$. Hence the
$D$-centered version suffices for those arguments.
\end{definition}

\begingroup
{The first growth lemma \cite[Lemma~2.1]{Safonov} provides a PDE-based mechanism for verifying (UC) for subsolutions that vanish on a prescribed boundary portion.}
\par\endgroup

\begin{remark}[Essential inequalities and countable covers]\label{rem:essential}
For supersolutions without a continuity assumption, all lower bounds below are essential lower bounds on open balls. Finite chains propagate such bounds between balls. The metric space has a countable base, so every open cover admits a countable subcover. Thus uniform local essential bounds yield an essential bound on the covered region with a single null exceptional set.

Pointwise notation for a measurable function $f$ may also be interpreted outside a fixed null set as follows. For each member $B$ of a countable base, discard the null set on which $f$ lies outside the essential range of $f$ on $B$. Outside their union, every open neighborhood $U$ of the point satisfies $\essinf_U f\le f(x)\le\esssup_U f$.
\end{remark}

\begin{lemma}[Essential strong minimum principle]\label{lem:strongmin}
Assume (Comp) and (CD). Let $G\subset D$ be a connected open set and let ${u}\in\mathcal P(G)$ be nonnegative almost everywhere. Either ${u}=0$ almost everywhere in $G$, or every nonempty open $W\Subset G$ satisfies
\begin{equation}\label{eq:strongmin}
 \essinf_W {u}>0.
\end{equation}
In the second case, every compact subset of $G$ has an open neighborhood on which ${u}$ has a positive essential lower bound.
\end{lemma}
\begin{proof}
If ${u}$ is not zero almost everywhere, some superlevel set $\{{u}\ge t\}$, with $t>0$, has positive measure. The differentiation theorem for doubling measures gives a density point of this set in $G$. Choose a sufficiently small ball centered there, with its $A_{\mathrm{CD}}$-dilate in $G$, on which the density is at least $\vartheta$. By (CD), ${u}$ has a positive essential lower bound on the concentric half-ball. Fix a point $y\in G$ and join the initial lower-bound ball to $y$ by a path whose image $K$ is compactly contained in $G$. Choose a common radius $r>0$ so small that the first ball $B(z_0,r)$ lies in the initial lower-bound ball and $B(z,2A_{\mathrm{CD}}r)\subset G$ for every $z\in K$. By uniform continuity of the path, choose finitely many centers $z_0,\ldots,z_N\in K$, with $z_N=y$ and $d(z_i,z_{i+1})<r$. For consecutive centers,
\[
\begin{aligned}
B(z_i,r)&\subset B(z_{i+1},2r)\subset B(z_i,4r),\\
\frac{\mu(B(z_i,r))}{\mu(B(z_{i+1},2r))}
&\ge\frac{\mu(B(z_i,r))}{\mu(B(z_i,4r))}\\
&\ge C_\mu^{-2}>\vartheta.
\end{aligned}
\]
If $u\ge a_i>0$ almost everywhere on $B(z_i,r)$, the density hypothesis of (CD) therefore holds in $B(z_{i+1},2r)$ at level $a_i$. The inclusion $B(z_{i+1},2A_{\mathrm{CD}}r)\subset G$ permits its application and gives
\[
 \essinf_{B(z_{i+1},r)}u\ge a_i/M.
\]
Starting with the initial lower bound and iterating along the finite chain gives a positive essential lower bound near $y$. A finite subcover of any compact subset of $G$ gives the final assertion and \eqref{eq:strongmin}.
\end{proof}

\begin{lemma}\label{lem:firstgrowth}
For each fixed $b\in(0,1)$, there exists $\epsilon_b>0$, determined by $b$, $M$, $C_\mu$, and the dilation factors, with the following property. {Let $U\subset D$ be open, $y\in D$, and $R>0$. Put $B=B(y,R)$ and assume $B(y,2R)\subset U$. If $u\in\mathcal C(U)$ is nonnegative and continuous, $c\ge0$, and}
\[
 \mu(B\cap\{{u}>c\})\le\epsilon_b\mu(B),
\]
then, for ${w}=({u}-c)_+$,
\[
 \esssup_{B(y,R/2)}{w}\le b\,\esssup_B {w}.
\]
\end{lemma}
\begin{proof}
By the interior hypothesis in Definition~\ref{def:UC}, (UC) holds for ${w}=({u}-c)_+$ because ${u}$ is a nonnegative continuous member of $\mathcal C$ and $c\ge0$. Fix $L\ge1$ with $2^{-L}\le b$. For each $z\in B(y,R/2)$ put $s_j=2^{-j-2}R$ for $j=0,\ldots,L-1$. Then $B(z,2s_j)\subset B(y,R)$. By doubling and the bounded distance between $z$ and $y$, every one of these balls satisfies
\[
 \mu(B(z,2s_j))\ge {c_*}\mu(B(y,R)).
\]
{Here $c_*>0$ depends only on $C_\mu$ and $L$.} Take $\epsilon_b<{c_*}/M$. {Hence,} the density of $\{{w}>0\}$ is then at most $1/M$ in each $B(z,2s_j)$. Applying (UC) with $t\downarrow0$ at the radii $s_0,\ldots,s_{L-1}$ gives
\[
 {w}(z)\le\esssup_{B(z,s_{L-1})}{w}
 \le2^{-L}\esssup_{B(z,2s_0)}{w}
 \le b\esssup_{B(y,R)}{w}.
\]
{{\color{black}Since $u$ is continuous, $w=(u-c)_+$ is continuous, and the bound at $z$ holds pointwise.} Taking the essential supremum over $z\in B(y,R/2)$ proves the assertion.}
\end{proof}

{The covering arguments of \cite[Section~7]{AFT} and \cite[Theorem~4.7]{DFGL08} amplify a fixed critical-density threshold. For the local version, the high superlevel sets are open by continuity, and (CD) bounds the covering radii. The covering step uses \cite[Theorem~4.1]{AFT}.}

\begin{lemma}[Power decay from critical density]\label{lem:powerdecay}
There are structural constants $\chi\in(0,1)$, $\sigma>0$, $C_{\mathrm{pd}}\ge1$, and $C_*>1$ with the following property. {\color{black}Let $u$ be a \emph{continuous nonnegative member of $\mathcal C$} on $B(x,C_*R)$.} Then, writing $m=\inf_{B(x,R)}{u}$, one has, for every $T\ge1$ and $\epsilon>0$,
\begin{equation}\label{eq:powerdecay}
\frac{\mu(\{{u}>T(m+\epsilon)\}\cap B(x,\chi R))}{\mu(B(x,R))}
 \le C_{\mathrm{pd}} T^{-\sigma}.
\end{equation}
The constants depend only on the measure-doubling constant, $M$, and the fixed ball dilations.
\end{lemma}
\begin{proof}
Choose $C_*=2+8A_{\mathrm{CD}}$. {In every application of (CD) below, the testing ball has center $y\in B(x,R)$ and radius at most $8R$. Indeed, the largest radius occurs when a lower bound on $B(y,b)$, with $b<2R$, is propagated using $B(y,4b)$. Hence the required testing dilations lie in $B(x,(1+8A_{\mathrm{CD}})R)\subset B(x,C_*R)$.} Let $g$ be a continuous nonnegative comparison supersolution on $B(x,C_*R)$ with $\inf_{B(x,R)}g=1$, and let $F\subset\{g>t\}\cap B(x,R)$ be nonempty and open, where $t>M$. For $y\in F$, the relative measure of $F$ in $B(y,4R)$ is less than $\vartheta$. Otherwise (CD) would give $g\ge t/M>1$ on $B(y,2R)\supset B(x,R)$, a contradiction.

To select a radius without a sphere-nullity assumption, put $r_k=4R\,2^{-k}$, $k\ge0$. Since $F$ is open, $B(y,r_k)\subset F$ for all sufficiently large $k$. Take the first $k=k(y)\ge1$ for which
\[
 \frac{\mu(F\cap B(y,r_k))}{\mu(B(y,r_k))}\ge\vartheta.
\]
Set $s_y=r_k$ and $b_y=2s_y$. {By the minimal choice of $k$,} the preceding dyadic ball has density strictly less than $\vartheta$, so
\begin{equation}\label{eq:stoppingdensity}
 0<b_y\le4R,\qquad
 \frac{\mu(F\cap B(y,s_y))}{\mu(B(y,s_y))}\ge\vartheta,\qquad
 \frac{\mu(F\cap B(y,b_y))}{\mu(B(y,b_y))}<\vartheta.
\end{equation}
By (CD) on $B(y,s_y)$, one has $g\ge t/M$ on $B(y,s_y/2)$. Since
\[
 \frac{\mu(B(y,s_y/2))}{\mu(B(y,2s_y))}\ge C_\mu^{-2}>\vartheta,
\]
a second application of (CD) on $B(y,2s_y)$ gives $g\ge t/M^2$ on $B(y,s_y)$. A third application on $B(y,4s_y)$, using $\mu(B(y,s_y))/\mu(B(y,4s_y))\ge C_\mu^{-2}>\vartheta$, gives the following bound on $B(y,b_y)$:
\begin{equation}\label{eq:stoppinglower}
 g\ge t/M^3.
\end{equation}
{Since $g$ is continuous,} these essential bounds hold pointwise. All testing dilations are contained in $B(x,C_*R)$.

Put $W=\bigcup_{y\in F}B(y,b_y)$. A countable subfamily covers $W$, since the ambient metric space is separable. Apply the $5$-covering lemma to this subfamily, obtaining disjoint balls $B_k$ whose fivefold dilates cover $W$. {By doubling,} $\mu(W)\le {C_\mu^3}\sum_k\mu(B_k)$. By the strict upper density bound in~\eqref{eq:stoppingdensity} and the disjointness of the selected balls,
\[
 \mu(W\setminus F)\ge\sum_k\mu(B_k\setminus F)
 \ge(1-\vartheta)\sum_k\mu(B_k)
 \ge\frac{1-\vartheta}{{C_\mu^3}}\mu(W).
\]
All sets have finite measure, being contained in $B(x,5R)$. Retain the countable subfamily $B(y_j,b_j)$ covering $W$, not only the disjoint selection. {Consequently,} with $q=1-(1-\vartheta)/{C_\mu^3}<1$ we obtain
\begin{equation}\label{eq:coverdata}
\begin{aligned}
 F&\subset\bigcup_jB(y_j,b_j),\\
 \mu(F\cap B(y_j,b_j))&<\vartheta\mu(B(y_j,b_j)),\\
 \mu(F)&\le q\mu\left(\bigcup_jB(y_j,b_j)\right).
\end{aligned}
\end{equation}
Moreover, $g\ge t/M^3$ on every covering ball, by~\eqref{eq:stoppinglower}.

If $t>M^3$, then $b_j<2R$: otherwise $B(x,R)\subset B(y_j,b_j)$ contradicts $\inf_{B(x,R)}g=1$. Starting from $B(y_j,b_j)$, repeatedly apply (CD) on four times the current radius. The lower bound loses a factor $M$ each time the known ball doubles. Stop at the first radius at least $2R$, which is less than $4R$. The resulting ball contains $B(x,R)$, and all testing dilations are contained in $B(x,C_*R)$. With $\gamma=\log_2M>0$, the number of doublings is $N=\lceil\log_2(2R/b_j)\rceil\le1+\log_2(2R/b_j)$, so
\[
 1\ge tM^{-(N+3)}\ge (t/M^4)(b_j/(2R))^\gamma.
\]
Consequently, since $M^{1/\gamma}=2$, {for every $t>M^3$, every covering radius satisfies}
\begin{equation}\label{eq:coverradius}
 b_j\le {32}R t^{-1/\gamma}.
\end{equation}

Fix $L>M^3$ and choose $P>M^3$ so large that
\[
 \sum_{n=0}^{\infty} {32}(PL^{n+1})^{-1/\gamma}<{\frac{1}{2}}.
\]
Define $R_0=R$ and
\[
 R_{n+1}=R_n-{32}R(PL^{n+1})^{-1/\gamma},\qquad
 F_n=B(x,R_n)\cap\{g>PL^n\}.
\]
{By the choice of $P$,} $R_n>R/2$. Apply \eqref{eq:coverdata} to $F_{n+1}$; the empty-set case requires no argument. By \eqref{eq:coverradius}, each covering ball lies in $B(x,R_n)$. The three critical-density applications give $g\ge PL^{n+1}/M^3>PL^n$ on that ball. Thus their union is contained in $F_n$, and
\[
 \mu(F_{n+1})\le q\mu(F_n).
\]
{Hence, induction} yields
\[
 \mu(B(x,R/2)\cap\{g>PL^n\})\le q^n\mu(B(x,R)).
\]
For $T\ge P$, choose $n$ with $PL^n\le T<PL^{n+1}$. With $\sigma=-\log q/\log L$, this proves a bound $C_{\mathrm{pd}}T^{-\sigma}$ for the normalized distribution of $g$ on $B(x,R/2)$. {The same bound holds} for $1\le T<P$ after enlarging $C_{\mathrm{pd}}$.

Finally, take $g=({u}+\epsilon)/(m+\epsilon)$. By (Comp), this is a continuous nonnegative comparison supersolution with infimum one on $B(x,R)$; it need not be a solution. {Therefore,} the inclusion $\{{u}>T(m+\epsilon)\}\subset\{g>T\}$ proves \eqref{eq:powerdecay} with $\chi=1/2$.
\end{proof}

\begin{lemma}[Arbitrary-density propagation]\label{lem:amplify}
{\color{black}Let $C_*$ and $\chi$ be the constants in Lemma~\ref{lem:powerdecay}.}
For each $\eta\in(0,1)$, there exist $\beta(\eta)>0$ and fixed ball-dilation factors such that,
for every nonnegative continuous member $u$ of $\mathcal C$ on $B(x, C_{*} R)$,
\begin{equation}\label{eq:amplification}
\mu(\{{u}\ge t\}\cap B(x,\chi R))\ge\eta\mu(B(x,\chi R))
\quad\Longrightarrow\quad
\inf_{B(x,R)}{u}\ge\beta(\eta)t.
\end{equation}
\end{lemma}
\begin{proof}
Let $m=\inf_{B(x,R)}{u}$. For $0<\epsilon<t$ and $0<t' < t$, apply \eqref{eq:powerdecay} with $T=t'/(m+\epsilon)$ if $T\ge1$; otherwise $m+\epsilon>t'$ already gives the conclusion. By doubling, $\mu(B(x,\chi R))/\mu(B(x,R))\ge c_\chi>0$, and hence
\[
\eta c_\chi\le C_{\mathrm{pd}}\left(\frac{m+\epsilon}{t'}\right)^\sigma.
\]
Letting $t'\uparrow t$ and $\epsilon\downarrow0$ proves the result with
\[
\beta(\eta)=\min\left\{1,\left(\frac{\eta c_\chi}{C_{\mathrm{pd}}}\right)^{1/\sigma}\right\}.
\]
\end{proof}

\begin{lemma}[Propagation from critical density]\label{lem:propagation}
{Let $U\subset D$ be open and let $u\in\mathcal P(U)$ be nonnegative. There exist $c,\gamma>0$, depending only on $C_\mu$, $M$, and $A_{\mathrm{CD}}$, such that for every $y\in U$, $R>0$, $t>0$, and $0<\rho\le1$,}
\begin{equation}\label{eq:propagation}
 {u}\ge t\text{ on }B(y,\rho R),\quad
 B(y,C_{\mathrm{prop}}R)\subset U
 \quad\Longrightarrow\quad
 \essinf_{B(y,3R)}{u}\ge c\rho^\gamma t.
\end{equation}
{Here $C_{\mathrm{prop}}=12A_{\mathrm{CD}}$.}
\end{lemma}
\begin{proof}
For $\rho$ bounded below by a fixed positive constant, a fixed number of applications of (CD) proves the estimate. For smaller $\rho$, apply (CD) on $B(y,4s)$ whenever ${u}\ge t_j$ on $B(y,s)$, and use
$\mu(B(y,s))/\mu(B(y,4s))\ge C_\mu^{-2}>\vartheta$.
{Hence,} ${u}\ge t_j/M$ on $B(y,2s)$. Starting with $s=\rho R$ and repeating at most $N\le C(1+\log(1/\rho))$ times, then applying finitely many comparable-ball changes to reach $B(y,3R)$, yields $M^{-N}t\ge c\rho^\gamma t$. Stop at the first doubled radius at least $3R$, which is less than $6R$. Every ball used in (CD) then has radius less than $12R$, and its required dilation is contained in $B(y,C_{\mathrm{prop}}R)$.
\end{proof}

\begin{proposition}[Harnack inequality from a weighted maximum]\label{prop:interior}
Under the ambient hypotheses, (Comp), (CD), and the interior form of (UC), {there exist structural constants $H\ge1$ and $C_{\mathrm{int}}>3+\max\{C_*/\chi,C_{\mathrm{prop}}\}$, depending only on $C_\mu$, $M$, and the prescribed analytic dilations, such that for every open $U\subset D$, every positive continuous $u\in\mathcal C(U)$, every $x\in X$, and every $R>0$ satisfying $B(x,C_{\mathrm{int}}R)\subset U$, the following scale-invariant Harnack inequality holds:}
\begin{equation}\label{eq:intH}
 \esssup_{B(x,R)}{u}\le H\essinf_{B(x,R)}{u}.
\end{equation}
\end{proposition}
\begin{proof}
Normalize $R=1$. {Choose} $\gamma$ from the lower-bound propagation estimate. Set
\[
 K=\max_{z\in{\overline{B(x,2)}}}(2-d(z,x))^\gamma {u}(z)
 =(2-d(z_0,x))^\gamma {u}(z_0),\qquad
 \rho={\frac{1}{2}}(2-d(z_0,x)).
\]
The weight vanishes on ${\overline{B(x,2)}\setminus B(x,2)}$, and ${u}$ is bounded and continuous on ${\overline{B(x,2)}}$. Thus a maximum is attained at some $z_0\in{B(x,2)}$. Also $\sup_{{B(x,1)}}{u}\le K$. On $B(z_0,\rho)$,
\[
 {u}\le\rho^{-\gamma}K=2^\gamma {u}(z_0).
\]
Since ${u}$ is a nonnegative continuous member of $\mathcal C$, the interior version of (UC) applies directly to $({u}-{u}(z_0)/2)_+$. The contrapositive of Lemma~\ref{lem:firstgrowth}, with $b<2^{-\gamma-1}$, gives a structural $\eta_0>0$ such that
\[
 \mu\bigl(B(z_0,\rho)\cap\{{u}\ge {u}(z_0)/2\}\bigr)
 \ge\eta_0\mu(B(z_0,\rho)).
\]
Apply Lemma~\ref{lem:amplify} with $R'=\rho/\chi$ and $\eta=\eta_0$. Its density ball is exactly $B(z_0,\rho)$, and its equation ball is $B(z_0,C_*\rho/\chi)$. {Since $d(x,z_0)<2$ and $\rho\le1$, the choice $C_{\mathrm{int}}>3+\max\{C_*/\chi,C_{\mathrm{prop}}\}$ places both this ball and the propagation balls inside $B(x,C_{\mathrm{int}})$.} Restriction of its conclusion to $B(z_0,{\rho/4})$ gives
\[
 \essinf_{B(z_0,{\rho/4})}{u}\ge\beta(\eta_0){u}(z_0)/2
 =\beta(\eta_0)(2\rho)^{-\gamma}K/2.
\]
 Apply \eqref{eq:propagation} with center $z_0$, outer radius $1$, and radius fraction ${\rho/4}\le1/4$. Since ${B(x,1)}\subset B(z_0,3)$, this propagates the bound into ${B(x,1)}$; the resulting factor $\rho^\gamma$ cancels the factor $\rho^{-\gamma}$. Thus $\essinf_{{B(x,1)}}{u}\ge cK$, while $\esssup_{{B(x,1)}}{u}\le K$. Undo the normalization to prove \eqref{eq:intH}.
\end{proof}

\begin{remark}[Nonnegative comparison functions]\label{rem:nonnegative}
The Harnack inequality~\eqref{eq:intH} also holds for a continuous nonnegative member $u$ of $\mathcal C$. Define the weighted maximum $K$ as in the proof of Proposition~\ref{prop:interior}. If $K=0$, then $u=0$ on $B(x,1)$, and the conclusion is immediate. If $K>0$, the maximizing point $z_0$ satisfies $u(z_0)>0$, and the proof of Proposition~\ref{prop:interior} applies. By (Comp), differences of positive solutions are continuous members of $\mathcal C$.
\end{remark}

\begin{lemma}[Interior oscillation from the Harnack inequality]\label{lem:intosc}
There exist $A_i>2$, $C_i\ge1$, and $\alpha_i\in(0,1)$, depending only on the interior Harnack data, such that if $u,v>0$ are solutions on $B(x,A_i T)\subset D$, then, for every $y,z\in B(x,T/4)$,
\begin{equation}\label{eq:intosc}
 \left|\frac{u(y)}{v(y)}-\frac{u(z)}{v(z)}\right|
 \le C_i\left(\frac{d(y,z)}{T}\right)^{\alpha_i}
 \osc_{B(x,T)}\frac uv.
\end{equation}
\end{lemma}
\begin{proof}
Write $f=u/v$ and let $H$ and $C_{\mathrm{int}}$ be the constants in \eqref{eq:intH}. Set $\theta=(2C_{\mathrm{int}})^{-1}$ and $\lambda=1-(2H^2)^{-1}$. Consider any ball $B(p,t)\subset B(x,T)$. Let $m=\inf_{B(p,t)}f$ and $L=\sup_{B(p,t)}f$. {These numbers are finite because $f$ is continuous on the compact set $\overline{B(x,T)}$.}

The functions ${w}_1=u-mv$ and ${w}_2=Lv-u$ are nonnegative continuous members of $\mathcal C$ on $B(p,t)$. Their Harnack inequalities follow from Remark~\ref{rem:nonnegative}. At the fixed point $p$, at least one satisfies ${w}_i(p)/v(p)\ge(L-m)/2$. Applying the Harnack inequality to ${w}_i$ and to $v$ on $B(p,\theta t)$ gives, for every $z\in B(p,\theta t)$,
\[
 \frac{{w}_i(z)}{v(z)}\ge H^{-2}\frac{{w}_i(p)}{v(p)}
 \ge\frac{L-m}{2H^2},
\]
because $B(p,C_{\mathrm{int}}\theta t)\subset B(p,t)$. Therefore
\begin{equation}\label{eq:intcontraction}
 \osc_{B(p,\theta t)}f\le\lambda\osc_{B(p,t)}f.
\end{equation}

Fix $y,z\in B(x,T/4)$. {The case $d(y,z)=0$ is immediate.} The ball $B(y,T/2)$ is contained in $B(x,T)$. If ${d(y,z)}\ge\theta T/2$, the original oscillation bound proves \eqref{eq:intosc} with a structural constant. Otherwise choose $n\ge1$ such that
\[
 \theta^{n+1}T/2\le {d(y,z)}<\theta^nT/2.
\]
Both points belong to $B(y,\theta^nT/2)$. Set $\alpha_*=\log\lambda/\log\theta>0$. {Therefore,} iterate \eqref{eq:intcontraction} about $y$ to obtain
\[
 |f(y)-f(z)|\le\lambda^n\osc_{B(x,T)}f
 \le\left(\frac{2{d(y,z)}}{\theta T}\right)^{\alpha_*}
 \osc_{B(x,T)}f.
\]
{Finally,} choose $\alpha_i=\min\{\alpha_*,1/2\}$. Since ${d(y,z)}/T<1$, decreasing the exponent preserves the estimate. This proves \eqref{eq:intosc}. The enlargement $A_i>2$ guarantees the neighborhood of $\overline{B(x,T)}$ used at the start.
\end{proof}

\section{Boundary estimates from geometry and {supersolution families}}\label{sec:boundary}

\subsection{Boundary hypotheses and main results}

{Fix $\Lambda_0=4\max\{\Lambda,C_{\mathrm{int}}\}$. The testing balls in Lemma~\ref{lem:chains} are chosen with dilation factor $C_{\mathrm{int}}$, and their dilates are contained in $B(\xi,\Lambda r)\subset B(\xi,\Lambda_0r)$. For the boundary results below, the center and scale conditions are stated in each theorem, lemma, and proposition.}

\begin{theorem}[One-sided boundary comparison on a diameter-uniform domain]\label{thm:bhp}
Assume the ambient setting, diameter uniformity, and (Comp), (CD), and (UC) on interior sets, and fix $\xi\in X\setminus D$ and $0<r<r_*/(4\Lambda_0)$ with $d(\xi,D)<r/16$. Let $u\ge0$ be a supersolution and $v>0$ a solution in $D\cap B(\xi,\Lambda_0r)$. Put $V=B(\xi,\Lambda_0r)$ and assume that (UC), with the common constant $M$, holds on $V$ for $(t u-v)_-$ for every $t\ge0$, on balls centered in $D$. Choose a fixed $c_a>0$ with $c_a<\kappa/(16\max\{A_{\mathrm{CD}},C_{\mathrm{int}}\})$, set $a=a_{\xi,r}$, and {define}
\[
 \mathcal Q_a(u):=\sup\bigl\{t\ge0:\mu({B(a,c_ar)}\cap\{u\ge t\})\ge\vartheta\mu({B(a,c_ar)})\bigr\}.
\]
Then, for almost every $x\in D\cap B(\xi,r)$,
\begin{equation}\label{eq:BHP}
 \mathcal Q_a(u)\le C\,v(a)\,\frac{u(x)}{v(x)}.
\end{equation}
{\color{black}Here $C$ is a structural constant.} If $u$ is a positive solution, $\mathcal Q_a(u)\asymp u(a)$ by \eqref{eq:intH}, and if (UC) also holds on $V$ for $(t v-u)_-$ for every $t\ge0$ on balls centered in $D$, then the usual two-sided, scale-invariant BHP follows. The quotient regularity is stated in Theorem~\ref{thm:quotient}.
\end{theorem}

\begin{theorem}[Scale-invariant two-point H\"older estimate]\label{thm:quotient}
Under the assumptions of Theorem~\ref{thm:bhp}, suppose $u$ and $v$ are positive solutions in $D\cap B(\xi,\Lambda_0r)$ and (UC) also holds on $V=B(\xi,\Lambda_0r)$ for $(t v-u)_-$ for every $t\ge0$ on balls centered in $D$. There are $\alpha\in(0,1)$ and $C\ge1$, independent of $\xi$, $r$, $u$, and $v$, such that
\begin{equation}\label{eq:twoholder}
 \left|\frac{u(x)}{v(x)}-\frac{u(y)}{v(y)}\right|
 \le C\left(\frac{d(x,y)}{r}\right)^\alpha
 \frac{u(a_{\xi,r})}{v(a_{\xi,r})},
\end{equation}
for every $x,y\in D\cap B(\xi,r/16)$. In particular, $u/v$ has a unique positive $\alpha$-H\"older extension to the metric closure of $D\cap B(\xi,r/16)$.
\end{theorem}

\subsection{Carleson estimate}\label{sec:comparison}
\begin{lemma}[Growth along a localized ball chain]\label{lem:pathgrowth}
{Let $\xi\in X\setminus D$ and $0<r<r_*/(4\Lambda_0)$ satisfy $d(\xi,D)<r/16$. Let $v>0$ be a continuous solution in $D\cap B(\xi,\Lambda_0r)$. Then there exist structural constants $C,\beta>0$ such that, for every $x\in D\cap B(\xi,4r)$,}
\begin{equation}\label{eq:pathgrowth}
 v(x)\le C\left(1+\frac r{\del(x)}\right)^\beta v(a_{\xi,r}).
\end{equation}
\end{lemma}
\begin{proof}
Apply Lemma~\ref{lem:chains}(i) with the dilation required by the interior Harnack inequality. Consecutive center values, and the endpoints, are compared in the concentric Harnack subballs constructed in Lemma~\ref{lem:packingchain}. The factor for each step is fixed; there are at most $C(1+\log^+(r/\del(x)))$ steps. All testing dilates lie in the equation neighborhood by the choice of $\Lambda_0$. {Hence,} exponentiating this step bound proves~\eqref{eq:pathgrowth}.
\end{proof}

\begin{lemma}[Power tails imply a local supremum bound]\label{lem:tailsup}
Let ${f}\ge0$ be defined on $D\cap V$ and satisfy (UC), on balls centered in $D$, on an ambient open neighborhood $V$ of $\overline{B(x_0,2R)}$. Suppose that, for every $T\ge T_0$,
\begin{equation}\label{eq:tail}
 \frac{\mu(D\cap\{{f}>T\}\cap B(x_0,2R))}{\mu(B(x_0,2R))}
 \le C_{\mathrm{tail}}(T_0/T)^\sigma,
\end{equation}
with $C_{\mathrm{tail}},\sigma,T_0>0$. Then $\esssup_{D\cap B(x_0,R)}{f}\le C T_0$, where $C$ depends only on the doubling and analytic data and $C_{\mathrm{tail}},\sigma$.
\end{lemma}
\begin{proof}
Since $X$ is proper and $\overline{B(x_0,2R)}\subset V$, a finite collection of balls centered in $D$, with their doubled radii contained in $V$, covers $D\cap B(x_0,2R)$. The qualitative boundedness in (UC) therefore gives $\esssup_{D\cap B(x_0,2R)}{f}<\infty$. Set $R_0=R$ and
\[
 R_{j+1}-R_j=\frac{6R}{\pi^2(j+1)^2},\qquad
 s_j=\frac{R_{j+1}-R_j}{8},\qquad
 S_j=\esssup_{D\cap B(x_0,R_j)}{f}.
\]
Then $R_j\uparrow2R$. {For every $y\in D\cap B(x_0,R_j)$, $B(y,4s_j)\subset B(x_0,R_{j+1})\subset V$. Hence (UC) applies at $y$ with radius $s_j$, and its outer essential supremum is at most $S_{j+1}$.} Doubling gives
\[
 \mu(B(y,2s_j))\ge c(s_j/R)^Q\mu(B(x_0,2R)).
\]
By \eqref{eq:tail}, the threshold $T_j=C T_0(R/s_j)^{Q/\sigma}$ makes the relative density of $\{{f}>T_j\}$ in each $B(y,2s_j)$ at most $1/M$. Apply (UC) on these balls and cover $D\cap B(x_0,R_j)$ by a countable subfamily of the sets $D\cap B(y,s_j)$. {Thus,} taking essential suprema gives $S_j\le MT_j+\frac12 S_{j+1}$. {Iterating from $j=0$ to $j=n$ gives}
\[
 {S_0\le M\sum_{j=0}^n2^{-j}T_j+2^{-n-1}S_{n+1}.}
\]
{Since $T_j/T_0\le C'(j+1)^{2Q/\sigma}$, the geometric series converges.} {Moreover,} the remainder vanishes by the previously assumed qualitative boundedness on the outer ball. {Therefore, $S_0\le CT_0$, as required.}
\end{proof}

\begin{proposition}[{Scale-invariant Carleson estimate}]\label{prop:CE}
{Assume the ambient hypotheses, diameter uniformity, and (Comp), (CD), and (UC) on interior sets. Let $\xi\in X\setminus D$ and $0<r<r_*/(4\Lambda_0)$ satisfy $d(\xi,D)<r/16$. Let $v>0$ be a solution in $D\cap B(\xi,\Lambda_0r)$.} Put $V=B(\xi,\Lambda_0r)$ and suppose that (UC) holds on $V$ for $v$, on balls centered in $D$. Then, for $a=a_{\xi,r}$,
\begin{equation}\label{eq:CE}
 \esssup_{D\cap B(\xi,2r)}v\le C_{\rm CE}v(a),
\end{equation}
{\color{black}where $C_{\rm CE}$ is a structural constant.}
\end{proposition}
\begin{proof}
Choose $0<\eta<\alpha_{\mathrm{lay}}$ from Lemma~\ref{lem:layer}, and let $\beta$ be the exponent in~\eqref{eq:pathgrowth}. Integrating the power $\eta/\beta$ of that estimate and using~\eqref{eq:localnegativepower} gives, by Chebyshev's inequality, for every $T\ge1$,
\[
 \frac{\mu(D\cap\{v>T v(a)\}\cap B(\xi,4r))}{\mu(B(\xi,4r))}
 \le C T^{-\eta/\beta}.
\]
{Therefore, Lemma~\ref{lem:tailsup} applied to $f=v$ on $D\cap V$ with $R=2r$ yields~\eqref{eq:CE}.}
\end{proof}

\subsection{Almost positivity and boundary comparison}\label{sec:bhp}

\begin{lemma}[Almost positivity from uniform geometry and density]\label{lem:almost}
{There exist structural constants $\varepsilon>0$, $a_0\in(0,1)$, and $N_0>1$ with the following property. Let $\xi\in X\setminus D$ and $0<r<r_*/(4\Lambda_0)$ satisfy $d(\xi,D)<r/16$. Put} $R_j=r(1+2^{-j})$, $t_j=\varepsilon r2^{-j}$, $D_j=D\cap B(\xi,R_j)$, and $F_j=D_j\cap\{\del\ge t_j\}$.
If $w\in\mathcal P(D\cap B(\xi,2r))$, (UC) holds on $B(\xi,2r)$ for $w_-$ on balls centered in $D$, and
$w\ge-1$ on $D_0$ and $w\ge N_0$ on $F_0$, {\color{black}then $w\ge0$ almost everywhere on $D\cap B(\xi,r)$.}
\end{lemma}
\begin{proof}
Let $H_{\mathrm{sc}}\ge1$ be the uniform lower-bound loss along the fixed short chains in Lemma~\ref{lem:chains}(ii), applied to the nonnegative comparison supersolution $w+1$. Fix an integer $L_0\ge1$ and set $a_0=2^{-L_0}<1/(2H_{\mathrm{sc}})$. Then choose $N_0$ large enough that $H_{\mathrm{sc}}^{-1}(N_0+1)-1\ge N_0a_0$. Let $C_{\mathrm{short}}$ be a fixed constant large enough to bound the distances from the initial center and the fixed interior-ball enlargements of those short chains by $C_{\mathrm{short}}t$. Choose $c_0>0$ such that
\begin{equation}\label{eq:locconstants}
 c_0<1/16,\qquad (\Lambda+8)c_0<1/8.
\end{equation}
Then choose $\varepsilon>0$ so small that
\begin{equation}\label{eq:epsilonchoice}
 {\varepsilon<\min\!\left\{\frac{1}{16(1+C_{\mathrm{short}})},\frac{1}{16},\frac{\kappa c_0}{16\,2^{L_0}}\right\}},\qquad
 C_\mu^2C_L\left(\frac{\varepsilon2^{L_0}}{c_0}\right)^{\alpha_{\mathrm{lay}}}\le\frac1M.
\end{equation}
All choices are structural and independent of $r$ and $j$.

{Suppose} $W\in\mathcal P(D_j)$ satisfies $W\ge-1$ on $D_j$, $W\ge N_0$ on $F_j$, and (UC) holds on $B(\xi,R_j)$ for $W_-$ on balls centered in $D$. On the region $F_{j+1}\cap\{\del\ge t_j\}$, we have $W\ge N_0$. If $x\in F_{j+1}$ and $t_j/2\le\del(x)<t_j$, Lemma~\ref{lem:chains}(ii) with $t=t_j$ connects $x$ to a ball in $F_j$ by at most a fixed number of interior balls. Because
\[
 d(x,\xi)<R_{j+1},\qquad (1+C_{\mathrm{short}})t_j
 <{\frac{1}{8}}r2^{-j}<R_j-R_{j+1}={\frac{1}{2}}r2^{-j},
\]
all these balls and their fixed Harnack dilates lie in $D\cap B(\xi,R_j)$. Applying (CD) along the chain to $W+1$ gives the uniform lower bound $H_{\mathrm{sc}}^{-1}(N_0+1)$ on an interior neighborhood of $x$, in the essential sense. A countable subcover and Remark~\ref{rem:essential} therefore give $W(x)+1\ge H_{\mathrm{sc}}^{-1}(N_0+1)$ for almost every such $x$. Hence $W\ge N_0a_0$ on $F_{j+1}$.

\begingroup
The negative part $W_-$ satisfies $0\le W_-\le1$ on $D_j$ and vanishes wherever $\del\ge t_j$ in $D_j$. Fix $x\in D_{j+1}$ with $\del(x)<t_j$ and choose $\zeta_x\in X\setminus D$ with $d(x,\zeta_x)=\del(x)$, using properness. For this fixed $j$, define $s_j=c_0 2^{-L_0}r2^{-j}$ and, for $0\le\ell<L_0$, set $\rho_\ell=2^\ell s_j$.
The choice \eqref{eq:epsilonchoice} ensures $d(x,\zeta_x)=\del(x)\le t_j\le s_j/16\le \rho_\ell/16$. {Moreover, $\zeta_x\in X\setminus D$, $d(\zeta_x,D)<\rho_\ell/16$, and $0<\rho_\ell\le(c_0/2)r2^{-j}<r_*$. The required neighborhood inclusion is}
\[
 B(\zeta_x,\Lambda \rho_\ell)\subset
 B\!\left(\xi,R_{j+1}+(\varepsilon+\Lambda c_0/2)r2^{-j}\right)
 \subset B(\xi,R_j).
\]
The same holds for $B(\zeta_x,4\rho_\ell)$ and $B(x,4\rho_\ell)$. Thus the boundary-layer estimate at $(\zeta_x,\rho_\ell)$ and the (UC) applications remain within the neighborhood where $W\ge-1$ and (UC) holds for ${W_-}$ on $D_j$.

In $D\cap B(x,2\rho_\ell)$, $\{{W_-}>0\}$ is contained in $\{\del<t_j\}$. The inclusion $B(x,2\rho_\ell)\subset B(\zeta_x,4\rho_\ell)$ and doubling, using $B(\zeta_x,4\rho_\ell)\subset B(x,5\rho_\ell)$, give
\begin{equation}\label{eq:densneg}
\begin{aligned}
&\frac{\mu(D\cap\{{W_-}>0\}\cap B(x,2\rho_\ell))}
      {\mu(B(x,2\rho_\ell))}\\
&\quad\le C_\mu^2
 \frac{\mu(D\cap B(\zeta_x,4\rho_\ell)\cap\{\del<t_j\})}
      {\mu(B(\zeta_x,4\rho_\ell))}\\
&\quad\le C_\mu^2 C_L\left(\frac{t_j}{\rho_\ell}\right)^{\alpha_{\mathrm{lay}}}\\
&\quad\le C_\mu^2 C_L
 \left(\frac{\varepsilon 2^{L_0}}{c_0}\right)^{\alpha_{\mathrm{lay}}}
 \le M^{-1}.
\end{aligned}
\end{equation}
Apply (UC) with $t\downarrow0$ successively at the radii $\rho_0,\rho_1,\ldots,\rho_{L_0-1}$. Since $B(x,2^{L_0}s_j)\subset B(\xi,R_j)$ and ${W_-}\le1$ there,
\[
 {W_-}(x)\le\esssup_{D\cap B(x,s_j)}{W_-}\le 2^{-L_0}\esssup_{D\cap B(x,2^{L_0}s_j)}{W_-}\le a_0
\]
for almost every $x\in D_{j+1}$, using the representative convention in Remark~\ref{rem:essential}. If $\del(x)\ge t_j$, this holds because ${W_-}(x)=0$. Therefore $W\ge-a_0$ on $D_{j+1}$ and $W\ge N_0a_0$ on $F_{j+1}$. The (UC) hypothesis passes to restrictions and positive rescalings, since $(\lambda W)_-=\lambda W_-$ for $\lambda>0$, and \eqref{eq:UD} is invariant under positive rescaling with the same $M$.

\endgroup

Start with $W_0=w$ and set $W_{j+1}=a_0^{-1}W_j$ after the one-step assertion. Induction shows $W_j\ge-1$ on $D_j$ and $W_j\ge N_0$ on $F_j$. Consequently $w\ge-a_0^j$ on $D_j$ for every $j$, and $D\cap B(\xi,r)\subset D_j$. Letting $j\to\infty$ proves $w\ge0$ there, up to sets of measure zero.
\end{proof}

\begin{proof}[Proof of Theorem~\ref{thm:bhp}]
Put $B_a=B(a,c_ar)$. Since $u$ is finite almost everywhere on $B_a$ and $0<\mu(B_a)<\infty$, continuity from above gives the following limit as $t\to\infty$:
\[
 \mu(B_a\cap\{u\ge t\})\longrightarrow0.
\]
Thus $\mathcal Q_a(u)<\infty$. If $\mathcal Q_a(u)=0$, inequality~\eqref{eq:BHP} follows from $u\ge0$. Suppose $\mathcal Q_a(u)>0$. For every $0<t<\mathcal Q_a(u)$, the definition of the supremum and monotonicity of the superlevel sets give
\[
 \mu(B_a\cap\{u\ge t\})\ge\vartheta\mu(B_a).
\]
Since $B(a,A_{\mathrm{CD}}c_ar)\subset D\cap B(\xi,\Lambda_0r)$ by the choice of $c_a$, (CD) implies
\[
 \essinf_{B(a,c_ar/2)}u\ge t/M.
\]
Letting $t\uparrow\mathcal Q_a(u)$ gives $\essinf_{B(a,c_ar/2)}u\ge\mathcal Q_a(u)/M$. Lemma~\ref{lem:chains}(i) propagates this to
$u(x)\ge c_\varepsilon\mathcal Q_a(u)$ on $\{x\in D\cap B(\xi,2r):\del(x)\ge\varepsilon r\}$, with $c_\varepsilon$ structural. For each point in the thick region, the corresponding chain includes a ball centered at that point; these balls therefore cover the region. A countable subcover makes the lower bound valid outside one null set, as in Remark~\ref{rem:essential}. With $C_{\rm CE}$ from Proposition~\ref{prop:CE}, form
\[
 w=K\frac{u}{\mathcal Q_a(u)}-\frac{v}{C_{\rm CE}v(a)}.
\]
By (Comp) and positive rescaling, $w\in\mathcal P$; $w\ge-1$ on $D\cap B(\xi,2r)$, and $0\le w_-\le v/(C_{\rm CE}v(a))$. The assumed (UC) for $(t u-v)_-$, including $t=0$ for the Carleson bound on $v$, applies to $w_-$ by positive rescaling on the balls centered in $D$. Choose $K$ from the thick-set lower bound so that $w\ge N_0$ there, with the $N_0,\varepsilon$ fixed by Lemma~\ref{lem:almost}. {Therefore,} that lemma gives $w\ge0$ on $D\cap B(\xi,r)$, proving \eqref{eq:BHP}. If (UC) also holds for $(tv-u)_-$ for every $t\ge0$, interchanging $u$ and $v$ in the same argument gives the two-sided BHP.
\end{proof}

\begin{corollary}[One-sided comparison for a fixed pair]\label{cor:fixedpair}
Keep the setting and notation of Theorem~\ref{thm:bhp}. There is a structural
constant $K_*\ge1$ such that, when ${\mathcal Q_a(u)>0}$, the conclusion
\eqref{eq:BHP} holds if the (UC) estimate on $V$, with the same $M$ and
balls centered in $D$, is assumed only for $v$ and
\[
 (t_*u-v)_-,\qquad t_* = \frac{K_* C_{\rm CE}v(a)}{{\mathcal Q_a(u)}},
\]
where $C_{\rm CE}$ is the constant of Proposition~\ref{prop:CE}.
For ${\mathcal Q_a(u)=0}$, \eqref{eq:BHP} is immediate.
The same assertion holds for $u\ge0$ in $\mathcal P$, provided the particular comparison $t_*u-v$ belongs to $\mathcal P$.
\end{corollary}
\begin{proof}
Fix $K_*$ as the structural factor required for the thick-set lower bound
in the proof of Theorem~\ref{thm:bhp}. That proof uses boundary (UC)
only for $v$ in the Carleson estimate and for
$w_-=(C_{\rm CE}v(a))^{-1}(t_*u-v)_-$ in Lemma~\ref{lem:almost}.
{Positive-rescaling invariance of (UC) transfers the assumed estimate for $(t_*u-v)_-$ to $w_-$. Thus the proof of Theorem~\ref{thm:bhp} applies with the stated hypotheses.}
\end{proof}

\subsection{Boundary H\"older continuity of quotients}\label{sec:quotient}

\begin{lemma}[Boundary oscillation, uniform in the center]\label{lem:boundosc}
Assume the hypotheses of Theorem~\ref{thm:bhp} for two positive solutions $u,v$. Suppose also that (UC) holds on $B(\xi,\Lambda_0r)$ for $(tv-u)_-$ for every $t\ge0$. Put
\[
 N=\frac{u(a_{\xi,r})}{v(a_{\xi,r})},\qquad
 S=\frac{r}{64\Lambda_0},\qquad {\rho_{\mathrm{bd}}}=(4\Lambda_0)^{-1},\qquad {L_{\mathrm{bd}}}=32/{\rho_{\mathrm{bd}}}.
\]
There are ${C_{\mathrm{bd}}}\ge1$ and ${\alpha_{\mathrm{bd}}}\in(0,1)$, depending only on the structural data, such that, for every $\zeta\in(X\setminus D)\cap B(\xi,r/8)$ and $0<s\le S$ with $s\ge {L_{\mathrm{bd}}} d(\zeta,D)$,
\begin{equation}\label{eq:boundosc}
 \osc_{D\cap B(\zeta,s)}\frac uv
 \le {C_{\mathrm{bd}}}\left(\frac{s}{S}\right)^{{\alpha_{\mathrm{bd}}}}N.
\end{equation}
\end{lemma}
\begin{proof}
Theorem~\ref{thm:bhp}, applied in both directions at $(\xi,r)$, gives the following estimate on $D\cap B(\xi,r)$:
\begin{equation}\label{eq:globalratio}
 C^{-1}N\le u/v\le CN.
\end{equation}
Let $\zeta\in(X\setminus D)\cap B(\xi,r/8)$. The initial quotient oscillation on $D\cap B(\zeta,S)$ is at most $CN$ by~\eqref{eq:globalratio}. All equation enlargements at scales at most $S$ lie in $B(\xi,\Lambda_0r)$.

Fix ${L_{\mathrm{bd}}} {d(\zeta,D)}\le t\le S$ with $t>0$, and write
\[
 m_t=\inf_{D\cap B(\zeta,t)}u/v,\qquad
 L_t=\sup_{D\cap B(\zeta,t)}u/v.
\]
{Since $d(\zeta,D)\le t/L_{\mathrm{bd}}<\rho_{\mathrm{bd}}t/16$ and $0<\rho_{\mathrm{bd}}t\le\rho_{\mathrm{bd}}S=r/(256\Lambda_0^2)<r_*/(4\Lambda_0)$, the center and radius satisfy the conditions of Corollary~\ref{cor:fixedpair}.} In $D\cap B(\zeta,t)$ the two functions
${w}_1=u-m_tv$ and ${w}_2=L_tv-u$ are nonnegative continuous members of $\mathcal C$ on the comparison neighborhood. At ${a_t=a_{\zeta,{\rho_{\mathrm{bd}}}t}}$,
\[
 \frac{{w}_1({a_t})}{v({a_t})}+\frac{{w}_2({a_t})}{v({a_t})}=L_t-m_t.
\]
If $L_t=m_t$, the contraction is immediate. Otherwise, choose $i\in\{1,2\}$ so that
\[
 \frac{w_i(a_t)}{v(a_t)}\ge\frac{L_t-m_t}{2}>0.
\]
Put $s=\rho_{\mathrm{bd}}t$. The anchor ball $B(a_t,c_as)$ and its $C_{\mathrm{int}}$-dilate lie in $D\cap B(\zeta,t)$ by the choices of $c_a$ and $\rho_{\mathrm{bd}}$. On $B(a_t,c_as)$, {\color{black}the Harnack inequality for nonnegative comparison functions in Remark~\ref{rem:nonnegative}} gives
\[
 H^{-1}w_i(a_t)\le w_i\le Hw_i(a_t).
\]
With $q_i=\mathcal Q_{a_t}(w_i)$ at the scale $(\zeta,s)$, the definition of the quantile now yields
\[
 H^{-1}w_i(a_t)\le q_i\le Hw_i(a_t).
\]
In particular, $q_i>0$. The equation enlargement lies in $B(\zeta,t)$ because $\Lambda_0s=t/4$. Apply Corollary~\ref{cor:fixedpair} to $w_i$ and $v$, with its prescribed coefficient
\[
 T_i=\frac{K_*C_{\rm CE}v(a_t)}{q_i}.
\]
Taking $t=0$ in the (UC) hypothesis for $(tu-v)_-$ gives (UC) for $v$.
For $i=1$, the required negative part is
\[
 (T_1{w}_1-v)_-=(1+T_1m_t)
 \left(\frac{T_1}{1+T_1m_t}u-v\right)_-,
\]
and the assumed (UC) estimate for $(tu-v)_-$ applies. For $i=2$, the
Carleson estimate for $v$ and $0\le {w}_2\le L_tv$ give
$q_2\le L_tC_{\rm CE}v({a_t})$. Thus $T_2L_t\ge K_*\ge1$, and
\[
 (T_2{w}_2-v)_-=T_2\left((L_t-T_2^{-1})v-u\right)_-.
\]
The assumed (UC) estimate for $(tv-u)_-$ applies since $L_t-T_2^{-1}\ge0$. These same identities, before taking negative parts, show that $T_i {w}_i-v\in\mathcal P$ by (Comp). Thus the extended form of Corollary~\ref{cor:fixedpair} applies, even though ${w}_i$ need not belong to $\Scone$ or $\Hcone$.
Consequently Corollary~\ref{cor:fixedpair}, followed by the interior
Harnack inequality and Remark~\ref{rem:nonnegative}, yields, for every $z\in D\cap B(\zeta,{\rho_{\mathrm{bd}}}t)$,
\[
 \frac{{w}_i(z)}{v(z)}\ge c\frac{{w}_i({a_t})}{v({a_t})}.
\]
{Therefore,} the uniform $c>0$ gives
\[
 L_{{\rho_{\mathrm{bd}}}t}-m_{{\rho_{\mathrm{bd}}}t}\le(1-c/2)(L_t-m_t).
\]
For a target ${L_{\mathrm{bd}}} {d(\zeta,D)}\le s\le S$, iterate at $t=S,{\rho_{\mathrm{bd}}}S,\ldots,{\rho_{\mathrm{bd}}}^{n-1}S$, where ${\rho_{\mathrm{bd}}}^{n+1}S<s\le{\rho_{\mathrm{bd}}}^nS$; no step is needed if $n=0$. {At each input scale $t\ge s$, $d(\zeta,D)\le s/L_{\mathrm{bd}}\le t/L_{\mathrm{bd}}<\rho_{\mathrm{bd}}t/16$, while $0<\rho_{\mathrm{bd}}t\le\rho_{\mathrm{bd}}S<r_*/(4\Lambda_0)$. Thus the center and radius conditions for Corollary~\ref{cor:fixedpair} hold at every step.} Monotonicity of oscillation gives~\eqref{eq:boundosc}, with ${\alpha_{\mathrm{bd}}}=\log(1-c/2)/\log{\rho_{\mathrm{bd}}}>0$, decreased if necessary. For ${d(\zeta,D)=0}$ the iteration has no lower cutoff. The (UC) hypotheses on the smaller neighborhoods are used only for the combinations verified above.
\end{proof}

\begin{proof}[Proof of Theorem~\ref{thm:quotient}]
Put
\[
f=\frac{u}{v},\qquad
N=\frac{u(a_{\xi,r})}{v(a_{\xi,r})},\qquad
\rho=d(x,y)>0.
\]
The case $x=y$ is immediate.
Let $S,L_{\mathrm{bd}}$ be as in Lemma~\ref{lem:boundosc}, and fix $K\ge8A_i+8$.
If $\rho\ge S/(L_{\mathrm{bd}}(K+1))$, the bound~\eqref{eq:globalratio} proves~\eqref{eq:twoholder}.
Hence suppose $\rho<S/(L_{\mathrm{bd}}(K+1))$.

Set $\delta=\del(x)$ and choose $\zeta_x\in X\setminus D$ with $d(x,\zeta_x)=\delta$, using properness. Since $\xi\in X\setminus D$ and $x\in B(\xi,r/16)$, we have $\delta<r/16$ and $\zeta_x\in B(\xi,r/8)$. Also $d(\zeta_x,D)\le\delta$. If $\delta\le K\rho$, both $x,y$ belong to $B(\zeta_x,{L_{\mathrm{bd}}}(K+1)\rho)$, where
\[
 {L_{\mathrm{bd}}} d(\zeta_x,D)\le {L_{\mathrm{bd}}}\delta<{L_{\mathrm{bd}}}(K+1)\rho<S.
\]
Lemma~\ref{lem:boundosc} therefore gives $|f(x)-f(y)|\le CN(\rho/S)^{{\alpha_{\mathrm{bd}}}}$.

If $\delta>K\rho$, put $T=\delta/(2A_i)$. Then $B(x,A_iT)\subset D\cap B(\xi,r)$ and $\rho<T/4$, so Lemma~\ref{lem:intosc} applies. If ${L_{\mathrm{bd}}}\delta\le S$, {the ball $B(x,T)$ lies in $B(\zeta_x,L_{\mathrm{bd}}\delta)$. Since $d(\zeta_x,D)\le\delta$ and $\zeta_x\in (X\setminus D)\cap B(\xi,r/8)$, Lemma~\ref{lem:boundosc} applies with radius $L_{\mathrm{bd}}\delta$. Hence}
\[
 |f(x)-f(y)|\le CN(\rho/\delta)^{\alpha_i}(\delta/S)^{{\alpha_{\mathrm{bd}}}}
 \le CN(\rho/S)^\alpha,\qquad \alpha=\min\{\alpha_i,{\alpha_{\mathrm{bd}}}\}.
\]
Here $\rho<\delta<S$, and the fixed factor ${L_{\mathrm{bd}}}^{{\alpha_{\mathrm{bd}}}}$ is absorbed in $C$. If ${L_{\mathrm{bd}}}\delta>S$, use~\eqref{eq:globalratio} for the interior oscillation to obtain
\[
 |f(x)-f(y)|\le CN(\rho/\delta)^{\alpha_i}\le CN(\rho/S)^\alpha,
\]
since $\delta>S/{L_{\mathrm{bd}}}$. As $S/r$ is structural, this proves~\eqref{eq:twoholder}. Uniform continuity gives the extension to the metric closure, and~\eqref{eq:globalratio} makes its values positive.
\end{proof}

\subsection{Qualitative comparison consequences}\label{sec:qualitative}
\begin{proposition}[Local weak minimum principle]\label{prop:weakmin}
Assume the ambient hypotheses, diameter uniformity, (Comp), and (CD). Let $V\subset X$ be open with $D\setminus V\ne\varnothing$, and put $U=D\cap V$. Suppose $w\in\mathcal P(U)$ and that (UC) holds on $V$ for $(w+c)_-$ for every $c\ge0$, for balls centered in $D$. Assume also that for every $\eta>0$ there is a compact set $K_\eta\Subset V$ such that, $\mu$-almost everywhere on $D\cap(V\setminus K_\eta)$,
\begin{equation}\label{eq:artificialboundary}
 w\ge-\eta.
\end{equation}
Then $w\ge0$ almost everywhere in $U$.
\end{proposition}
\begin{proof}
First suppose $w_-=0$ almost everywhere on $D\cap(V\setminus K)$ for some compact $K\Subset V$. By (UC), $w_-$ is essentially bounded on $D$ in every ambient ball whose required enlargement lies in $V$.
Since $K$ is compact, the set $K \cap D$ is totally bounded.
Choose finitely many such balls covering this set, using total boundedness and the positive distance from $K$ to $X\setminus V$. Together with the vanishing of $w_-$ on $D\cap(V\setminus K)$, this gives $m:=\esssup_U w_-<\infty$. Suppose $m>0$ and set ${u}=w+m\ge0$. Each component of $U$ meets $D\cap(V\setminus K)$: a path in $D$ to a point of $D\setminus V$ leaves $K$ before leaving $V$. On this nonempty open intersection, ${u}\ge m$. Lemma~\ref{lem:strongmin} and a finite cover therefore give a positive essential lower bound for ${u}$ near every compact subset of $U$.

For each $\xi\in K\cap\partial D$, choose $0<s_\xi<r_*$ with $\overline{B(\xi,4s_\xi)}\subset V$. Fix $0<\tau<1$ so small that $C_\mu^4 C_L\tau^{\alpha_{\mathrm{lay}}}\le M^{-1}$ and put $t_\xi=\tau s_\xi$. {The set} ${\overline{B(\xi,2s_\xi)}\cap\{\del\ge t_\xi\}}$ {is compactly contained in $U$. Hence ${u}\ge c_\xi>0$ essentially near this set; if the set is empty, choose any $c_\xi>0$.} Choose $0<b_\xi<\min\{m,c_\xi\}$ and put ${w}_\xi=({u}-b_\xi)_-$ on $U$. Since $m-b_\xi\ge0$, (Comp) makes ${u}-b_\xi=w+(m-b_\xi)$ a comparison supersolution. Moreover, $({u}-b_\xi)_-=(w+m-b_\xi)_-$, so the assumed (UC) estimate applies directly to ${w}_\xi$ on balls centered in $D$. We have ${w}_\xi\le b_\xi$, and its positive set in $D\cap B(\xi,2s_\xi)$ lies, up to a null set, in $\{\del<t_\xi\}$.

For any $x\in D\cap B(\xi,s_\xi)$ put $\rho=s_\xi/4$. Since $B(x,2\rho)\subset B(\xi,2s_\xi)$, Lemma~\ref{lem:layer} and doubling give
\[
\begin{aligned}
&\frac{\mu(D\cap\{{w}_\xi>0\}\cap B(x,2\rho))}
      {\mu(B(x,2\rho))}\\
&\quad\le C_\mu^4
 \frac{\mu(D\cap B(\xi,4s_\xi)\cap\{\del<t_\xi\})}
      {\mu(B(\xi,4s_\xi))}\\
&\quad\le C_\mu^4 C_L\tau^{\alpha_{\mathrm{lay}}}\le M^{-1}.
\end{aligned}
\]
The comparison of denominators uses $B(\xi,4s_\xi)\subset B(x,5s_\xi)$ and $2\rho=s_\xi/2$. Since $B(x,4\rho)\subset B(\xi,2s_\xi)\subset V$, the $D$-centered (UC), with its level tending to zero, implies
\[
 \esssup_{D\cap B(x,\rho)}{w}_\xi
 \le{\frac{1}{2}}\esssup_{D\cap B(x,2\rho)}{w}_\xi\le b_\xi/2.
\]
The balls $D\cap B(x,\rho)$ cover $D\cap B(\xi,s_\xi)$. {Hence,} taking a countable subcover and combining the essential estimates shows ${u}\ge b_\xi/2$ essentially on $D\cap B(\xi,s_\xi)$.

Since $K\cap\partial D$ is compact, finitely many of the neighborhoods $B(\xi,s_\xi)$ cover it. Let $O$ be their union. The remaining set
\[
 (K\cap\overline D)\setminus O\Subset U
\]
is compact and has positive distance from $\partial D$ and from $X\setminus V$. Lemma~\ref{lem:strongmin} gives a positive essential lower bound for ${u}$ on a neighborhood of this compact remainder. Combining it with the finitely many bounds obtained near $\partial D$ yields ${u}\ge c>0$ almost everywhere on $K\cap D$. Outside $K$, ${u}\ge m$. Consequently $w_-\le m-\min\{c,m/2\}<m$ almost everywhere in $U$, contradicting the definition of $m$.

To treat the general condition~\eqref{eq:artificialboundary}, apply the first part to $w+\eta$. By (Comp), this is a comparison supersolution; the (UC) hypothesis for its nonnegative shifts follows from the hypothesis for $w$, and its negative part vanishes outside $K_\eta$. Thus $w\ge-\eta$ almost everywhere. Letting $\eta\downarrow0$ proves the assertion.
\end{proof}

\begin{remark}[Condition near $\partial V$]\label{rem:weakminboundary}
For $V\Subset X$, condition~\eqref{eq:artificialboundary} is equivalent to
\[
 \lim_{\rho\downarrow0}\esssup_{D\cap V\cap\{d(x,X\setminus V)<\rho\}}w_-=0,
\]
where the supremum over the empty set is zero. Indeed, compact subsets of $V$ have positive distance from $X\setminus V$, whereas the closed inner cores of $\overline V$ are compactly contained in $V$. The hypothesis $D\setminus V\ne\varnothing$ ensures that each connected component of $D\cap V$ meets $D\cap(V\setminus K)$ for every compact $K\Subset V$, as used in the proof of Proposition~\ref{prop:weakmin}.
\end{remark}

\begin{proposition}[Local comparison with a solution]\label{prop:localcomparison}
In the setting of Proposition~\ref{prop:weakmin}, let $u\in\Scone(U)$ and $v\in\Hcone(U)$. If (UC) holds on $V$ for $(u-v+c)_-$ for every $c\ge0$ and $u-v$ satisfies \eqref{eq:artificialboundary}, then $u\ge v$ almost everywhere in $U$.
\end{proposition}
\begin{proof}
Condition (Comp) gives $u-v\in\mathcal P(U)$. {Therefore,} Proposition~\ref{prop:weakmin} applies to this difference and gives $u\ge v$ almost everywhere in $U$.
\end{proof}

\section{Examples covered by the framework}\label{sec:examples}

{ Unless otherwise specified, in the linear examples let $\Scone(U)$ and $\Hcone(U)$ be the supersolutions and solutions, respectively, of the indicated equation on $U$, and set $\mathcal P(U)=\Scone(U)$. For equations with no zeroth-order term or with a nonnegative zeroth-order term in the sign convention $-\Delta+c$, linearity and the supersolution inequality for nonnegative constants give (Comp).}

\subsection{Elementary examples}\label{sec:examples-elementary}

\begin{example}[Concavity and averaging on intervals]\label{ex:concave}
\begingroup
On $X=\mathbb R$ with Lebesgue measure, let $D=(0,1)$ and let $\Scone(U)$ consist of finite functions concave on each component of $U$; then $\Hcone(U)$ consists of componentwise affine functions.
{ For a nonnegative concave function on $B(x,4s)$, concavity bounds its infimum on $B(x,s/2)$ from below by a fixed multiple of its supremum on $B(x,s)$. This gives (CD). The comparison functions in $\mathcal C(U)$ are affine on each component of $U$. Hence their positive truncations are convex on interior intervals, and the local upper estimate for convex functions gives interior (UC) by Section~\ref{sec:A2upper}.

Let $u$ be nonnegative and concave, and let $v$ be positive and affine with zero boundary value at an endpoint of $D$. For $t\ge0$, the function $(v-tu)_+$ is nonnegative and convex, is bounded above by $v$, and tends to zero at that endpoint. Its extension by zero across the endpoint is convex. Thus the same local upper estimate gives boundary (UC). Since the interval is diameter-uniform, Theorem~\ref{thm:bhp} applies. If both solutions vanish at the same endpoint, their quotient is constant.
}\endgroup
\end{example}

\subsection{Divergence-form equations and Dirichlet forms}\label{sec:examples-divergence}

\begin{example}[Uniformly elliptic divergence-form equations]\label{ex:elliptic}
\begingroup
{ Let $X=\mathbb R^n$ with Euclidean distance and Lebesgue measure, and let $D\subsetneq X$ be diameter-uniform. Consider}
$L=-\operatorname{div}(\bm{A}\nabla\cdot)$, where $\bm{A}$ is bounded,
{ measurable, and satisfies
\[
 \bm A(x)\zeta\cdot\zeta\ge\lambda_{\mathrm{ell}}|\zeta|^2,
 \qquad |\bm A(x)|\le\Lambda_{\mathrm{ell}}
\]
for almost every $x\in D$ and every $\zeta\in\mathbb R^n$, where $0<\lambda_{\mathrm{ell}}\le\Lambda_{\mathrm{ell}}$. The matrix need not be symmetric.}
The weak supersolution cone satisfies (Comp); the classical weak
Harnack inequality and local subsolution estimates give (CD) and interior (UC),
respectively \cite[Chapter~8]{GT01}. Let $V$ be an open boundary
neighborhood. Suppose that a positive comparison solution $v$ satisfies
$\eta v\in H_0^1(D)$ for every nonnegative
$\eta\in C_c^\infty(V)$.
Lemma~\ref{lem:A2traceenergy} and Section~\ref{sec:A2div} then give
(UC) on $V$ for $(v-{tu})_+$, uniformly for ${t}\ge0$
and nonnegative weak supersolutions $u$. If both comparison functions
are positive solutions with the same trace condition, (UC) holds
{ in both directions. Hence Theorems~\ref{thm:bhp} and~\ref{thm:quotient} give the boundary comparison and quotient estimates, with constants depending only on $n$, the diameter-uniformity constant, and $\Lambda_{\mathrm{ell}}/\lambda_{\mathrm{ell}}$. The Laplacian is the case $\bm{A}=I$.}
For bounded Euclidean uniform domains, Corollary~\ref{cor:A2capacitary}
identifies the localized zero-trace assumption with a finite-energy
{ condition given by \eqref{eq:A2capcondition}.}
\endgroup
\end{example}

\begingroup
{ The following identities are related to Doob's $h$-transform \cite{Doob57}. For a positive $C^2$ function $h$, direct differentiation gives
\[
\begin{aligned}
 -\Delta(hz)+\frac{\Delta h}{h}hz&=-h^{-1}\operatorname{div}(h^2\nabla z),\\
 -\Delta(hz)+\operatorname{div}((\nabla\log h)hz)&=-\operatorname{div}(h\nabla z).
\end{aligned}
\]
{ In both cases multiplication by $h$ preserves the quotient of two solutions. The equations below are understood in the weak sense, with solutions in $H^1_{\mathrm{loc}}(D)$. We use the Euclidean metric and Lebesgue measure. For the interior singularities we apply the interior theory to the transformed divergence-form equation. For the boundary singularities we use the Schr\"odinger equation with nonnegative potential, including the equation obtained by transforming the drift equation in Example~\ref{ex:drift}(ii).}

}\begin{example}[Critical Schr\"odinger potentials]\label{ex:schrodinger}
\mbox{}\par\nobreak
\begin{enumerate}[label=(\roman*),itemsep=0.35\baselineskip]
\item \emph{Interior singularity.}
Let $n\ge2$, $D=B(0,1)\subset\mathbb R^n$, and choose
$h(x)=1+\varepsilon/\log(e/|x|)$ for $x\ne0$, with $h(0)=1$ and
$0<\varepsilon<1$. The resulting equation is
\[
 -\Delta u+
 \frac{\varepsilon\big((n-2)\log(e/|x|)+2\big)}
 {|x|^2\log^2(e/|x|)\big(\log(e/|x|)+\varepsilon\big)}\,u=0.
\]
The potential is nonnegative and belongs to $L^{n/2}(D)$, with its
norm tending to zero as $\varepsilon\downarrow0$. The positive weak solution $h$ is not H\"older continuous
at $0$: for every $\alpha>0$,
\begin{equation}\label{eq:examplenonholder}
 \frac{|h(x)-h(0)|}{|x|^\alpha}
 =\frac{\varepsilon}{|x|^\alpha\log(e/|x|)}\longrightarrow\infty
 \qquad(x\to0).
\end{equation}
{ The transformation is valid also at the origin. Indeed, writing $r=|x|$ and $\ell=\log(e/r)$, we have $|\nabla h|\le C/(r\ell^2)$ and $0\le\Delta h\le C/(r^2\ell^2)$. The flux through $\partial B(0,r)$ tends to zero as $r\downarrow0$, {\color{black}so the classical Laplacian of $h$ on $D\setminus\{0\}$ represents its distributional Laplacian on $D$.} For compactly supported Sobolev functions in $D$, the weight $r^{-2}\ell^{-2}$ is controlled by the Dirichlet energy: this follows from the Hardy inequality for $n>2$ and, for $n=2$, from the one-dimensional Hardy inequality after the radial change of variable $\ell=\log(e/r)$. By localization, multiplication by $h$ and by $h^{-1}$ preserves $H^1_{\mathrm{loc}}(D)$, and the potential term defines a locally bounded bilinear form on $H^1$. Testing the equation for $u$ with $h\varphi$ and the equation for $h$ with $u\varphi$ and subtracting gives
\[
 \int_D h^2\nabla(u/h)\cdot\nabla\varphi\,dx=0
\]
for every $\varphi\in C_c^\infty(D)$. Conversely, the same identities give the equation for $u$ from that for $u/h$. Since $1\le h\le1+\varepsilon$, the transformed equation is uniformly elliptic. Thus its positive solutions are locally H\"older continuous, by the interior estimates in Example~\ref{ex:elliptic}. The quotient of any two positive weak solutions of the original equation is therefore locally H\"older continuous.}

\item \emph{Flat boundary singularity.}
Let $D=B(0,1)\cap\{x_n>0\}\subset\mathbb R^n$, $n\ge2$,
and choose $h(x)=x_n^p$, $p>1$. The resulting equation is
\[
 -\Delta u+\frac{p(p-1)}{x_n^2}\,u=0.
\]
{ Let $V$ be an open neighborhood of a point of the flat boundary, avoiding the spherical boundary. For each positive comparison solution $v$ assume, for every compact $K\subset V$,
\[
 \int_{D\cap K}\left(|\nabla v|^2+v^2+x_n^{-2}v^2\right)\,dx<\infty.
\]
For the original Schr\"odinger equation the nonnegative potential gives (Comp). On interior balls with $B(x,4s)\subset D$, the rescaled potential is bounded on $B(x,2s)$ by a constant depending only on $p$. The weak Harnack inequality for uniformly elliptic equations with bounded zeroth-order coefficient gives (CD) \cite[Chapter~8]{GT01}. Nonnegative solutions and their positive truncations are subsolutions of $-\Delta$, giving interior (UC).

For a nonnegative supersolution $u$, the positive part $w=(v-tu)_+$, $t\ge0$, is a weak subsolution of $-\Delta$ in $D\cap V$: on its positive set the potential contributes with the favorable sign, and the positive-part inequality follows by truncation. Moreover, $0\le w\le v$. The energy assumption implies \eqref{eq:A2HV} locally on the flat boundary. The cutoff argument in Section~\ref{sec:A2div} yields \eqref{eq:A2zerotrace} and the subsolution inequality for the zero extension. Hence (UC) holds on $V$ uniformly in $t$. If both positive solutions satisfy the energy condition, Theorems~\ref{thm:bhp} and~\ref{thm:quotient} give the BHP and H\"older continuity of their quotient. No $A_2$ assumption on $x_n^{2p}$ is needed.}

\item \emph{Polar boundary point.}
Let $D=B(0,1)\setminus\{0\}\subset\mathbb R^n$, $n\ge2$,
and choose $h(x)=|x|^a$, $a>0$. The resulting equation is
\[
 -\Delta u+\frac{a(a+n-2)}{|x|^2}\,u=0.
\]
{ The capacity-density condition fails at $0$. Let $V$ be a neighborhood of $0$ compactly contained in $B(0,1)$. For each positive solution $v$ assume finite local Schr\"odinger energy, namely
\[
 \int_{D\cap K}\left(|\nabla v|^2+v^2+|x|^{-2}v^2\right)\,dx<\infty
\]
for every compact $K\subset V$. Both this punctured ball and the half-ball in (ii) are uniform domains. The interior estimates in (ii) apply with the bound on $\delta_D^2 a(a+n-2)/|x|^2$ in place of the bound for the flat boundary potential. For $w=(v-tu)_+$, the subsolution inequality for $-\Delta$ and $0\le w\le v$ again give \eqref{eq:A2zerotrace} by the cutoff argument of Section~\ref{sec:A2div}, now with $\delta_D(x)=|x|$ near $0$. Its extension across $0$ is therefore a weak subsolution of $-\Delta$ in $V$. This gives boundary (UC) and, for two positive solutions satisfying the energy condition, the BHP and H\"older regularity of the quotient.} The weight $|x|^{2a}$ lies outside $A_2$ when
$a\ge n/2$, but the example allows every $a>0$.
\end{enumerate}
\end{example}

{ For the boundary behavior of Schr\"odinger operators with scale-critical potentials on Lipschitz domains, see, e.g., \cite{An12,Mar19}. Example~\ref{ex:schrodinger}(iii) gives boundary comparison at a point where the capacity density condition fails. {\color{black}For potentials more singular than the inverse-square boundary-distance scale, even an isolated boundary point} may correspond to infinitely many minimal Martin boundary points; see \cite{Nakai75,Suzuki84}.}

\begin{example}[Critical drifts]\label{ex:drift}
\mbox{}\par\nobreak
\begin{enumerate}[label=(\roman*),itemsep=0.35\baselineskip]
\item \emph{Interior singularity.}
{ Let $D$ and $h$ be as in Example~\ref{ex:schrodinger}(i). The drift equation is}
\[
 -\Delta u+
 \operatorname{div}\!\left(
 \frac{\varepsilon x}
 {|x|^2\log(e/|x|)\big(\log(e/|x|)+\varepsilon\big)}\,u
 \right)=0.
\]
The drift belongs to $L^n(D)$, and its norm tends to zero as
$\varepsilon\downarrow0$.
The function $h$ is a positive weak solution, but it is not H\"older
continuous at $0$ by \eqref{eq:examplenonholder}. Nevertheless,
{ the quotient of any two positive weak solutions is locally H\"older continuous. Indeed, the Sobolev multiplier estimate in Example~\ref{ex:schrodinger}(i) yields $u/h\in H^1_{\mathrm{loc}}(D)$ and the identity $\nabla u-u\nabla\log h=h\nabla(u/h)$ almost everywhere. Thus the weak equation is precisely the uniformly elliptic equation $-\operatorname{div}(h\nabla(u/h))=0$.}

\item \emph{Flat boundary singularity.}
Let $D=B(0,1)\cap\{x_n>0\}\subset\mathbb R^n$, $n\ge2$,
and choose $h(x)=x_n^{-\beta}$, $\beta>0$. The resulting equation is
\[
 -\Delta u-\beta\,\partial_n\!\left(\frac{u}{x_n}\right)=0.
\]
With $v=x_n^{\beta/2}u$, this becomes
\[
 -\Delta v+\frac{\beta(\beta+2)}{4x_n^2}\,v=0,
\]
which is the Schr\"odinger equation in Example~\ref{ex:schrodinger}(ii)
{ with $p=1+\beta/2$. This identity holds locally for weak solutions since $x_n^{\beta/2}$ is smooth and positive in $D$. For each of the two positive solutions of the drift equation, require its transform to satisfy the local energy condition in Example~\ref{ex:schrodinger}(ii). The transformation preserves their quotient. Hence the BHP and H\"older estimate there apply for every $\beta>0$, without an $A_2$ assumption on $x_n^{-\beta}$. The constants may depend on $\beta$.}
\end{enumerate}
\end{example}
\endgroup

\begin{example}[A spectral fractional problem through its extension]\label{ex:spectral}
\begingroup
Let $\Omega\subset\mathbb R^n$ be a bounded uniform domain, let
$0<s<1$, and write $\mathcal A_\Omega=-\Delta_\Omega$ for the
Dirichlet Laplacian. Suppose that nonzero, nonnegative
$u,v\in\operatorname{Dom}(\mathcal A_\Omega^{s/2})$ satisfy
$(\mathcal A_\Omega^{s/2}w,\mathcal A_\Omega^{s/2}\varphi)_{L^2(\Omega)}=0$
for $w=u,v$ and every $\varphi\in C_c^\infty(\Omega\cap B(\xi,R))$,
where $\xi\in\partial\Omega$. \begingroup
For $w=u,v$, let $w_N$ denote the partial sums of its Dirichlet
eigenfunction expansion, and let $W_N$ be the corresponding Poisson
extensions \cite[Theorem~1.1]{ST10}. The spectral expansion gives
\[
 \int_{\Omega\times(0,\infty)}
 y^{1-2s}|\nabla W_N|^2\,dx\,dy
 =\frac{2^{1-2s}\Gamma(1-s)}{\Gamma(s)}
 \|\mathcal A_\Omega^{s/2}w_N\|_{L^2(\Omega)}^2.
\]
The same identity applied to differences shows that $W_N$ converges in
weighted gradient energy. Since $\Omega$ is bounded, the positive first
Dirichlet eigenvalue also controls the weighted $L^2$ norm, so the limit
belongs to $H^1(\Omega\times(0,\infty),y^{1-2s}\,dx\,dy)$ and has trace $w$
at $y=0$. The positivity-preserving Dirichlet heat semigroup and the Poisson
formula give $W\ge0$ for nonnegative $w$.

The weighted weak formulation for the extension is
\[
\begin{aligned}
 &\int_{\Omega\times(0,\infty)}
 y^{1-2s}\nabla W\cdot\nabla\Phi\,dx\,dy\\
 &\quad=\frac{2^{1-2s}\Gamma(1-s)}{\Gamma(s)}
 (\mathcal A_\Omega^{s/2}w,
 \mathcal A_\Omega^{s/2}\Phi(\cdot,0))_{L^2(\Omega)}
\end{aligned}
\]
for every smooth compactly supported test function $\Phi$ on
$\Omega\times[0,\infty)$. Thus the assumed weak equation makes the
conormal term vanish on $\Omega\cap B(\xi,R)$. For a test function meeting the base $y=0$, apply this identity on the upper half-cylinder and its reflection, then add the two identities. { Both conormal terms vanish by the assumed weak equation.} Hence the even extensions $U,V$ satisfy the weighted weak equation
\[
 -\operatorname{div}(|y|^{1-2s}\nabla W)=0,\qquad W=U,V,
\]
near $(\xi,0)$ in
$D=\Omega\times\bigl(-2\operatorname{diam}\Omega,2\operatorname{diam}\Omega\bigr)$.
{ For each finite spectral sum, the extension belongs to $H_0^1(\Omega)$ on every horizontal section. Multiplying by a smooth ambient cutoff $\eta$ and extending by zero across the lateral boundary gives a function in the weighted ambient energy space (approximate each of the finitely many spatial eigenfunctions in $H_0^1(\Omega)$ by compactly supported smooth functions).

The convergence in weighted gradient energy and weighted $L^2$ then gives}
\[
 (\eta U)^0,(\eta V)^0\in H_0^1(B((\xi,0),R_1),|y|^{1-2s}\,dx\,dy)
\]
whenever $\eta\in C_c^\infty(B((\xi,0),R_0))$, $0<R_0<R_1<R$, and the superscript $0$ denotes zero extension. In particular, cutoffs localize the zero-trace property without prescribing values on $\partial B((\xi,0),R)$.
{ Since $\Omega$ is connected, the Dirichlet heat semigroup is positivity improving. The Poisson formula therefore gives positive extensions for $y>0$. On the portion of $y=0$ where the even extensions satisfy the equation, weighted interior regularity gives continuous representatives and the strong minimum principle gives positivity \cite{FKS82}. Their traces there agree almost everywhere with $u$ and $v$, respectively.}
\endgroup

Equip $X=\mathbb R^{n+1}$ with the Euclidean metric and
$d\mu=|y|^{1-2s}\,dx\,dy$. The weight belongs to $A_2$, the
measure is doubling, and $D$ is uniform. The weighted weak
supersolution cone satisfies (Comp); weighted interior estimates
\cite{FKS82} give (CD) and interior (UC). The zero-trace property
and the weighted argument of Section~\ref{sec:A2div} give
(UC) near $\partial D$ for $(tU-V)_-$ and $(tV-U)_-$, uniformly in $t\ge0$. Hence Theorems~\ref{thm:bhp} and~\ref{thm:quotient},
restricted to $y=0$, yield the boundary Harnack principle and the
two-point H\"older estimate for $u/v$ near $\xi$, at scales
$r<c\min\{R,\operatorname{diam}\Omega\}$, where $c>0$ depends only on the structural dilation factors and the uniformity constant.
\endgroup
\end{example}

\begin{example}[Harmonic functions on the Sierpi\'{n}ski gasket]\label{ex:gasket}
Let $K$ be the standard Sierpi\'{n}ski gasket in $\mathbb R^2$, with outer
side length one, Euclidean distance $d(x,y)=|x-y|$, and canonical
self-similar measure $\mu_K$. Let $V_0=\{q_0,q_1,q_2\}$ be its outer
vertices, and put $D=K\setminus V_0$. The measure $\mu_K$ is doubling.
{ The uniformity of $D$ is also noted in \cite{BCHOTW26},
where parabolic Anderson models are studied on inner uniform domains with
Dirichlet boundary conditions and on uniform domains with Neumann boundary
conditions.}

\paragraph{Diameter uniformity.}
The intrinsic shortest-path metric $d_g$ on $K$ satisfies
\begin{equation}\label{eq:gasket-geodesic}
 |x-y|\le d_g(x,y)\le8|x-y|
\end{equation}
for all $x,y\in K$; see \cite[Proposition~8.10]{CIL08}.
Fix $x,y\in D$, and let $\gamma$ be a shortest path between them,
parametrized by arc length. Every subarc of $\gamma$ is also shortest.
In particular, for each $z\in\gamma$,
\begin{equation}\label{eq:gasket-subarc}
 \ell(\gamma[x,z])\le8|x-z|,\qquad
 \ell(\gamma[z,y])\le8|z-y|.
\end{equation}

For an outer vertex $q_i$, write $C_{i,m}=F_i^m(K)$ for the corner
cell of side length $2^{-m}$. The only junctions between $C_{i,m}$ and the rest of the gasket
are its two vertices other than $q_i$, which we denote by
$p_{i,m}$ and $p'_{i,m}$. The straight segment
$[p_{i,m},p'_{i,m}]$ belongs to $K$. If both endpoints of a shortest
path lie outside $C_{i,m}$, any portion of the path passing through
$C_{i,m}$ must connect these two junction vertices. Such a portion
has length at least $2^{-m}$, and equality forces it to be the
straight segment $[p_{i,m},p'_{i,m}]$. Returning to the same junction
vertex would contradict minimality. Consequently, the part of a
shortest path in $C_{i,m}$ with endpoints outside that cell lies on
its opposite edge.

If $\gamma$ passed through $q_i$, choose $m$ large enough that
$x,y\notin C_{i,m}$. The preceding observation would force its part
inside the cell to lie on the opposite edge, which does not contain
$q_i$. Thus $\gamma\subset D$.

{ Fix $z\in\gamma$, and choose $i$ such that $\delta_D(z)=|z-q_i|$. If $\delta_D(z)<\sqrt3/4$, choose a dyadic side length $a=2^{-m}$ with $2\delta_D(z)/\sqrt3<a\le4\delta_D(z)/\sqrt3$. The point $z$ is in the corner cell $C_{i,m}$ and not on its opposite edge, since that edge is at distance $\sqrt3a/2$ from $q_i$. Thus at least one endpoint of $\gamma$ lies in the cell. If it is $x$, then \eqref{eq:gasket-subarc} gives $\ell(\gamma[x,z])\le8a\le32\delta_D(z)/\sqrt3$; the estimate for the other endpoint is the same. If $\delta_D(z)\ge\sqrt3/4$, the shorter subarc has length at most $\ell(\gamma)/2\le4|x-y|\le4$. In both cases the shorter subarc has diameter at most $32\delta_D(z)$. Since $\operatorname{diam}\gamma\le8|x-y|$, condition~\eqref{eq:diamuniform} holds with $A=32$.

}\paragraph{Local mean-value inequality.}
Let $V$ be relatively open in $K$, and let $f$ be a nonnegative,
locally bounded, locally finite-energy subharmonic function on $V$.
We claim that there is a structural constant $C_0\ge1$ such that,
whenever $B_K(x,4r)\subset V$ and $0<r\le1/16$,
\begin{equation}\label{eq:gasket-local-mv}
 \sup_{B_K(x,r)} f\le C_0
 \left(\frac{1}{\mu_K(B_K(x,2r))}
 \int_{B_K(x,2r)}f^2\,d\mu_K\right)^{1/2}.
\end{equation}
Here $B_K$ denotes a Euclidean ball intersected with $K$.

For $i\in\{0,1,2\}$, consider the unbounded blow-up
\[
 K_\infty^{(i)}=\bigcup_{m\ge0}\bigl(q_i+2^m(K-q_i)\bigr)
\]
with its standard strongly local resistance form and its
self-similar measure $\mu_\infty$ normalized by
$\mu_\infty|_K=\mu_K$. The sub-Gaussian heat-kernel upper bound
and conservativeness on this unbounded gasket, together with
volume doubling, give $(FK)_\Psi$ and $(Gcap_{\le})_\Psi$ by
\cite[Theorem~1.3]{GHL15}, with
$\Psi(r)=r^{d_w}$ and $d_w=\log5/\log2$.
\cite[Theorem~6.3]{GHL15} consequently gives the squared
mean-value inequality for bounded nonnegative subharmonic functions
on balls of $K_\infty^{(i)}$.

Because $4r\le1/4$, the ball $B_K(x,4r)$ contains at most one
outer vertex. Choose this vertex as $q_i$ if it is present, and
otherwise choose any outer vertex. The other two outer vertices
are outside $B_K(x,4r)$. The only attachment points between
$K$ and $K_\infty^{(i)}\setminus K$ are these two vertices.
Choose a finite-energy cutoff $0\le\eta\le1$ on $K$ that equals
one on $\overline{B_K(x,2r)}$ and has support in $B_K(x,3r)$.
The function $\eta f$ has finite energy and is zero near both
attachment points. Extend it by zero from $K$ to
$K_\infty^{(i)}$, and denote the extension by $\widetilde f$.
It belongs to the global energy domain and agrees with $f$ on
$B_K(x,2r)$. { Write $B_\infty$ for balls in the unbounded gasket.} { To see that $B_\infty(x,2r)=B_K(x,2r)$, write $q_j,q_k$ for the two attachment vertices. The new cells in the blow-up are contained in the union of the two translated cones with vertices $q_j,q_k$ generated by $q_j-q_i$ and $q_k-q_i$. At either attachment vertex the angle between this cone and the directions into the original triangle is at least $\pi/3$. Thus the distance from $x\in K$ to the cone at $q_l$ is at least $(\sqrt3/2)|x-q_l|$ for $l=j,k$. As $|x-q_l|\ge4r$, neither cone meets $B(x,2r)$.}
{ Strong locality now shows that $\widetilde f$ is subharmonic on $B_\infty(x,2r)$. Applying the mean-value inequality of \cite[Theorem~6.3]{GHL15}}
to $\widetilde f$ on this ball proves~\eqref{eq:gasket-local-mv}.
The essential supremum supplied by that theorem agrees with
the supremum of the continuous resistance-form representative.

The same inequality, with a different structural constant, holds
for $0<r\le1/4$ whenever $B_K(x,4r)\subset V$. Indeed, if
$1/16<r\le1/4$, apply~\eqref{eq:gasket-local-mv} at each
$z\in B_K(x,r)$ with radius $\rho=r/16$.
The enlarged ball $B_K(z,4\rho)$ lies in $B_K(x,2r)\subset V$.
By doubling, $\mu_K(B_K(z,2\rho))$ is bounded below by a fixed
multiple of $\mu_K(B_K(x,2r))$. Taking the supremum over $z$
gives the asserted bound. For interior (UC) balls one has
$4r\le\delta_D(x)\le1$; for the balls used below to verify (UC) near $\partial D$, $V$ avoids $q_1,q_2$, which likewise gives
$4r\le1$. Thus these are all the radii needed here.

Let $(\mathcal E,\mathcal F)$ be the standard resistance form
\cite{Kig01}. Harmonicity on an open set $U\subset D$ means
$\mathcal E(h,\varphi)=0$ for test functions with compact support
in $U$; for functions of finite global energy that are harmonic
on all of $D$, one may take
\[
 \mathcal F_0=\{\varphi\in\mathcal F:\varphi|_{V_0}=0\}
\]
as the test space. Let $\Hcone(U)=\Scone(U)$ consist of the
continuous locally finite-energy harmonic functions on $U$,
with $\mathcal P(U)$ the corresponding harmonic comparison
family. Comparison differences satisfy (Comp), and the elliptic
Harnack inequality implies (CD) for nonnegative harmonic functions.
If $h\in\mathcal C(U)$ and $c\ge0$, the convex truncation
$(h-c)_+$ is subharmonic. The local mean-value estimate and the
level-set splitting in Section~\ref{sec:A2upper} therefore give
interior (UC).

For the boundary condition near $q_0$, choose a relatively open
neighborhood $V\subset K$ with $V\cap\{q_1,q_2\}=\varnothing$.
Let $u,v\in\mathcal F\cap C(K)$ be nonnegative functions harmonic
on $D$, with $u(q_0)=v(q_0)=0$, and let $t\ge0$. Set
$g=v-tu$ and $w=g_+$. Then $0\le w\le v$ and $w(q_0)=0$.
For every $\varepsilon>0$, the function
$w_\varepsilon=(w-\varepsilon)_+$ vanishes in a neighborhood
of $q_0$. On $D$, it is subharmonic by the convex-truncation
property; hence it is subharmonic on all of $V$, including $q_0$.
By the Markov property of the resistance form,
\[
 \|w_\varepsilon-w\|_{L^{\infty}(K)}\le\varepsilon,
 \qquad
 \mathcal E(w_\varepsilon,w_\varepsilon)
 \le\mathcal E(w,w)\le\mathcal E(g,g).
\]
For each cutoff $\eta\in\mathcal F$ supported compactly in $V$,
the functions $\eta w_\varepsilon$ have uniformly bounded energy.
Their weak limit in the energy space is $\eta w$, by uniform
convergence. Testing the subharmonic inequalities with nonnegative
functions supported in the region where $\eta=1$ and passing to
the weak limit proves that $w$ is subharmonic on $V$,
including at $q_0$.

{ Since $\mu_K(V_0)=0$, estimate~\eqref{eq:gasket-local-mv} for $w$ has the form of \eqref{eq:A2localupper} on the balls required in (UC). The boundedness in Definition~\ref{def:UC} follows from $w\le v$ and continuity of $v$ on $K$. The level-set splitting in Section~\ref{sec:A2upper} therefore proves (UC) with a constant independent of $u$, $v$, and $t$.
}Taking $t=0$ gives the corresponding estimate for $v$, and
interchanging $u$ and $v$ gives (UC) for $(tv-u)_-$.
Thus the abstract boundary comparison and quotient estimates
apply to this class of gasket harmonic functions near $q_0$.
For related Dirichlet-space boundary Harnack results, see
\cite{GS11,LSC14,BM19,Chen25}.
The harmonic extensions of the vertex data $(0,1,1)$ and
$(0,1,2)$ provide nontrivial examples.
\end{example}

\subsection{Nondivergence-form equations}\label{sec:examples-nondivergence}

Bauman \cite{Bauman84} established boundary comparison for positive solutions of linear nondivergence-form equations in Lipschitz domains.
{ For boundary comparison with unbounded drift on John and uniform domains, see \cite{Safonov,KS11,KS14}.}

\begin{example}[Linear nondivergence-form equations]\label{ex:nondivlinear}
Let $n\ge2$, let $X=\mathbb R^n$ with Lebesgue measure, let $D$ be diameter-uniform, and let $V$ be an open boundary neighborhood.
Consider the operator
\begin{equation}\label{eq:nondivoperator}
 Lf=\sum_{i,j=1}^n a^{ij}D_{ij}f+\sum_{i=1}^n b^iD_if.
\end{equation}
{\color{black}Assume that $a^{ij}:D\cap V\to\mathbb R$ are measurable and $a^{ij}=a^{ji}$. For some constant $0<\nu\le1$, assume that the following holds for almost every $x\in D\cap V$ and every $\zeta\in\mathbb R^n$:}
\[
 \nu|\zeta|^2\le\sum_{i,j=1}^n a^{ij}(x)\zeta_i\zeta_j
                  \le\nu^{-1}|\zeta|^2.
\]
Assume also
\[
 b\in L^n(D\cap V;\mathbb R^n).
\]
Let $u,v\in W^{2,n}_{\mathrm{loc}}(D\cap V)\cap C(D\cap V)$ be nonnegative functions satisfying $Lu\le0$ and $Lv=0$ almost everywhere in $D\cap V$.
{ Assume in addition that $v$ extends continuously to $\overline D\cap V$ and vanishes on $\partial D\cap V$.}

For open $U\subset D\cap V$, take $\Scone(U)=\mathcal P(U)$ to be the continuous strong supersolutions of $Lw\le0$ in $W^{2,n}_{\mathrm{loc}}(U)$, and let $\Hcone(U)$ be the strong solutions of $Lw=0$. Linearity and $L1=0$ give (Comp), and $\mathcal C=\Hcone$. The density-to-lower-bound argument in \cite[Lemma~2.6]{Safonov} gives (CD) with constants determined by $n$, $\nu$, and $\|b\|_{L^n(D\cap V)}$. Under the change of variables $z=x+sy$, for every measurable set $E\subset\mathbb R^n$ such that $x+sE\subset D\cap V$, the rescaled drift $\widetilde b(y)=s\,b(x+sy)$ satisfies
\[
 \|\widetilde b\|_{L^n(E)}
 =\|b\|_{L^n(x+sE)}
 \le\|b\|_{L^n(D\cap V)}.
\]
The interior and boundary forms of (UC) follow from the first growth lemma as detailed in Section~\ref{sec:A2nondiv}, uniformly for ${t}\ge0$.

Consequently, when $v>0$, Theorem~\ref{thm:bhp} gives one-sided boundary comparison with every nonnegative strong supersolution $u$. If both $u$ and $v$ are positive strong solutions vanishing continuously on $\partial D\cap V$, (UC) also holds for $(tv-u)_-$, and Theorem~\ref{thm:quotient} gives the scale-invariant two-point H\"older estimate for $u/v$. The constants depend only on the geometric data, $n$, $\nu$, and $\|b\|_{L^n(D\cap V)}$.

\end{example}

\begingroup
\begin{theorem}[Nonnegative potentials]\label{thm:nonnegativepotentials}
{ Let $D$, $V$, and $L$ be as in Example~\ref{ex:nondivlinear}, with ellipticity constant $\nu$ and $b\in L^n(D\cap V;\mathbb R^n)$. Assume $c\ge0$ and $\delta_D^2c\in L^\infty(D)$, and set $L_c=L-c$. Let $u,v\in W^{2,n}_{\mathrm{loc}}(D\cap V)\cap C(D\cap V)$ satisfy $u\ge0$, $v>0$, $L_cu\le0$, and $L_cv=0$ almost everywhere in $D\cap V$. Assume that $v$ extends continuously as zero on $\partial D\cap V$. Then the one-sided boundary comparison in Theorem~\ref{thm:bhp} holds. If $u$ is also a positive solution extending continuously as zero on $\partial D\cap V$, the two-sided boundary Harnack principle and the quotient estimate of Theorem~\ref{thm:quotient} hold. The constants depend only on the geometric data, $n$, $\nu$, $\|b\|_{L^n(D\cap V)}$, and $\|\delta_D^2c\|_{L^\infty(D)}$.}
Moreover, there exist constants $C_{\mathrm H}\ge1$ and $\gamma_{\mathrm H}>0$,
depending only on $n$, $\nu$, and $\|b\|_{L^n(D\cap V)}$, such that every positive strong
solution $w$ of $L_cw=0$ on $B(x,4R)\subset D\cap V$ satisfies
\begin{equation}\label{eq:nondivpotentialharnack}
 \sup_{B(x,R)}w\le C_{\mathrm H}
 \exp\!\left(\gamma_{\mathrm H}\sqrt{\|\delta_D^2c\|_{L^{\infty}(D)}}\right)
 \inf_{B(x,R)}w.
\end{equation}
\end{theorem}

\begin{proof}
For $a\ge0$, linearity and $c\ge0$ give (Comp):
\[
 L_c(w+a)=L_cw-ca\le0.
\]
We verify (CD), (UC), and the interior Harnack estimate in two stages.

\medskip\noindent
{ \emph{Stage 1. Lower-bound propagation for small lower-order terms.}
Write $B_j=B(0,j)$ in the normalized balls of this proof. After rescaling to $B_4$, assume
\[
 \|b\|_{L^n(B_4)}\le s_0,\qquad \|c\|_{L^\infty(B_4)}\le\kappa_c,
\]
where $s_0$ and $\kappa_c$ are sufficiently small constants depending only on $n$ and $\nu$.

}{ Suppose that $w\ge0$ and $L_cw\le0$ in $B_4$, with $w\ge1$ in $B_1$. On $B_4\setminus\overline{B_1}$, take
\[
 \psi(x)=\frac{|x|^{-k}-4^{-k}}{1-4^{-k}},
\]
where $k>0$ is chosen so that $\nu(k+1)\ge(n-1)\nu^{-1}$. The radial and tangential eigenvalues of $D^2\psi$ give $L_0\psi\ge0$, where $L_0$ is the second-order part of $L$. Moreover, $0\le\psi\le1$ and $\psi\le w$ on the two boundary spheres.} On
$E=\{\psi>w\}$,
\[
 L(\psi-w)\ge b\cdot\nabla\psi-cw
 \ge b\cdot\nabla\psi-\kappa_c,
 \qquad
 \sup_E(\psi-w)
 \le C(n,\nu)(s_0+\kappa_c).
\]
Taking the two thresholds smaller than the minimum of $\psi$
on a fixed intermediate annulus gives $w\ge\lambda_0$ on
$\partial B_{r_0}$ for some $r_0\in(2,4)$ and structural $\lambda_0>0$.
In the enclosed ball, set
$E_0=B_{r_0}\cap\{w<\lambda_0\}$. There
\[
 L(\lambda_0-w)\ge-cw\ge-\kappa_c\lambda_0,
 \qquad
 \sup_{E_0}(\lambda_0-w)
 \le C\kappa_c\lambda_0\le\frac{\lambda_0}{2}.
\]
Hence
\[
 \inf_{B_1}w\ge1
 \quad\Longrightarrow\quad
 \inf_{B_2}w\ge\frac{\lambda_0}{2}.
\]
{ This estimate is uniform under the small-coefficient assumptions.}

\medskip\noindent
\emph{Stage 2. Estimates for arbitrary lower-order terms.}

{ First, we prove the upper contraction for a bounded right-hand side. Let $G\subset B(x,2R)$ be open, and suppose $f\in W^{2,n}_{\mathrm{loc}}(G)\cap C(\overline G\cap B(x,2R))$ satisfies $0\le f\le1$, $Lf\ge-K_0R^{-2}$ almost everywhere in $G$, and $f=0$ on $\partial G\cap B(x,2R)$. Here $K_0\ge0$. Rescale to $x=0$ and $R=1$. Choose $N$ so large that $N^{-1/n}\|b\|_{L^n(D\cap V)}$ is smaller than a constant to be fixed below (the norm before rescaling is used for this choice). Partition $B_2\setminus\overline{B_1}$ into $N$ shells of width $1/N$. There is one shell $E$ such that the rescaled drift has $L^n(E)$ norm at most $N^{-1/n}\|b\|_{L^n(D\cap V)}$. In the following estimates $b$ denotes the rescaled drift. Let $r_1$ be the middle radius of this shell and put $\rho=1/(4N)$.

For each $y\in G\cap\partial B_{r_1}$, {\color{black}apply the Alexandrov--Bakelman--Pucci (ABP) estimate \cite[Theorem~1.1]{Safonov}} to $f(z)-|z-y|^2/\rho^2$ on $G\cap B(y,\rho)$. This function is nonpositive on the boundary, and
\[
 L\left(f(z)-\frac{|z-y|^2}{\rho^2}\right)
 \ge-K_0-2n\nu^{-1}\rho^{-2}-2\rho^{-1}|b(z)|.
\]
Since $B(y,\rho)\subset E$, the estimate gives
\[
 f(y)\le C\left((\rho K_0+\rho^{-1})|G|^{1/n}+\|b\|_{L^n(E)}\right).
\]
Here $C$ depends only on $n$ and $\nu$ when the shell drift norm is at most one. Choose $N$ so that the drift contribution is at most $1/4$. The ABP estimate on $G\cap B_{r_1}$, with the full drift bound, then gives
\[
 \sup_{G\cap B_1}f\le\frac14+C(1+K_0)|G|^{1/n},
\]
where $C$ may additionally depend on $\|b\|_{L^n(D\cap V)}$ and $N$. Taking $\varepsilon>0$ sufficiently small and returning to the original scale proves
\begin{equation}\label{eq:nondivpotentialfirstgrowth}
 |G|\le\varepsilon|B(x,2R)|
 \quad\Longrightarrow\quad
 \sup_{G\cap B(x,R)}f\le\frac12.
\end{equation}
The threshold depends only on $n$, $\nu$, $\|b\|_{L^n(D\cap V)}$, and $K_0$.

}{ {\color{black}For lower-bound propagation, rescale to $B_4$. Let $w$ be nonnegative and continuous in $B_4$, with $w\in W^{2,n}_{\mathrm{loc}}(B_4)$, $L_cw\le0$ almost everywhere, and $\inf_{B_1}w\ge1$. Assume first that $\|c\|_{L^\infty(B_4)}\le\kappa_c$, where $\kappa_c$ is the threshold of Stage~1, to be decreased below.} Choose an integer $m\ge2$ such that $m>(2\|b\|_{L^n(D\cap V)}/s_0)^n$. Choose $m$ equally spaced directions on a unit circle in a two-dimensional coordinate plane. Around the radial segments with these directions and radial coordinate between $1/2$ and $7/2$, take tubes of radius $c_n/m$, with $c_n>0$ small enough that the tubes are disjoint. Partition $B_3\setminus\overline{B_2}$ into $m$ disjoint shells of width $1/m$. One tube $T$ and one shell $E$ satisfy
\[
 \|b\|_{L^n(T)},\ \|b\|_{L^n(E)}\le m^{-1/n}\|b\|_{L^n(B_4)}<s_0/2.
\]
Write $r_1\in(2,3)$ for the middle radius of $E$. Starting on the axis of $T$ at radius $3/4$, follow this axis to $\partial B_{r_1}$, and then follow a shorter great-circle arc on $\partial B_{r_1}$ to any $y$ on that sphere. Choose $\rho_0=c'_n/m$, with $c'_n>0$ sufficiently small that the balls of radius $4\rho_0$ along the radial segment lie in $T$, and those along the spherical arc lie in $E$. The first ball of radius $\rho_0$ is contained in $B_1$. Choose successive centers $z_0,\ldots,z_J$ at arc length at most $\rho_0$, including the terminal point $z_J=y$, so that $J+1\le C(n)m$. For each center $z_j$, therefore,
\[
 \|\rho_0 b(z_j+\rho_0\cdot)\|_{L^n(B_4)}\le s_0,
 \qquad
 \rho_0^2\|c\|_{L^\infty(B(z_j,4\rho_0))}\le\kappa_c.
\]
}{ Since}
$B(z_{j+1},\rho_0)\subset B(z_j,2\rho_0)$,
iteration gives, uniformly for $y\in\partial B_{r_1}$,
\[
 { w(y)\ge(\lambda_0/2)^{J+1}}
 { \ge(\lambda_0/2)^{C(n)m}}
 =:{\lambda_1}>0.
\]
Here $m$ and ${\lambda_1}$ depend only on $n$, $\nu$,
and $\|b\|_{L^n(D\cap V)}$.
On $E_0=B_{r_1}\cap\{w<{\lambda_1}\}$,
\[
 L({\lambda_1}-w)\ge-cw\ge-\kappa_c{\lambda_1}.
\]
The ABP estimate, with the full drift norm, gives
\[
 \sup_{E_0}({\lambda_1}-w)
 \le C(n,\nu,\|b\|_{L^n(D\cap V)})\kappa_c{\lambda_1}.
\]
{\color{black}Decrease $\kappa_c>0$, if necessary, allowing it to depend also on $\|b\|_{L^n(D\cap V)}$,} so that the last expression is at most
${\lambda_1}/2$. Hence,
\[
 \inf_{B_1}w\ge1
 \quad\Longrightarrow\quad
 \inf_{B_2}w\ge\frac{{\lambda_1}}{2}.
\]

Now let the normalized potential bound $K_0$ be arbitrary. Suppose
$B(x,4r)\subset D\cap V$, $w\ge0$, $L_cw\le0$ there, and
$r^2\|c\|_{L^{\infty}(B(x,4r))}\le K_0$. Set
\[
 \eta=\min\left\{\frac1{16},
 \frac18\sqrt{\frac{\kappa_c}{1+K_0}}\right\},
 \qquad \rho=\eta r.
\]
For every testing ball $B(z,4\rho)\subset B(x,4r)$,
\[
 \begin{aligned}
 \rho^2\|c\|_{L^{\infty}(B(z,4\rho))}
 &\le\eta^2K_0\le\kappa_c,\\
 \|\rho b(z+\rho\cdot)\|_{L^n(B_4)}
 &=\|b\|_{L^n(B(z,4\rho))}
 \le\|b\|_{L^n(D\cap V)}.
 \end{aligned}
\]
A chain of at most $C(n)\eta^{-1}$ overlapping $\rho$-balls,
whose $4\rho$-dilates remain in $B(x,4r)$, gives
\begin{equation}\label{eq:nondivpotentialdoubleball}
 \inf_{B(x,r)}w\ge t
 \quad\Longrightarrow\quad
 \inf_{B(x,2r)}w\ge\lambda t,
 \qquad
 \lambda\ge\left(\frac{{\lambda_1}}{2}\right)^{C(n)/\eta}>0.
\end{equation}
The constant $\lambda$ depends only on $n$, $\nu$,
$\|b\|_{L^n(D\cap V)}$, and $K_0$.

{ We can now verify (CD). Suppose $w\ge0$ and $L_cw\le0$ in $B(x,2R)\subset D\cap V$, with $R^2\|c\|_{L^\infty(B(x,2R))}\le K_0$. For any ball $B(z,4s)\subset B(x,2R)$, apply \eqref{eq:nondivpotentialfirstgrowth} to $f=1-w/t$ on $G=B(z,2s)\cap\{w<t\}$, where $t>0$. On $G$,
\[
 Lf=-t^{-1}Lw\ge-t^{-1}cw\ge-K_0s^{-2}.
\]
The function $f$ vanishes on the level-set boundary inside $B(z,2s)$. Hence, if $|G|\le\varepsilon|B(z,2s)|$, then $w\ge t/2$ on $B(z,s)$. The double-ball estimate~\eqref{eq:nondivpotentialdoubleball} propagates this lower bound along overlapping balls whose fourfold dilates are contained in $B(x,2R)$; both estimates have constants controlled by $K_0$.

To obtain a lower bound from an arbitrary positive density, cover $B(x,R)$ by the cubes of a grid of side length $R/(100\sqrt n)$ that meet $B(x,R)$. Their union is contained in $B(x,5R/4)$ and their total volume is at most a dimensional multiple of $|B(x,R)|$. If $|\{w\ge t\}\cap B(x,R)|\ge\theta|B(x,R)|$, one such cube $Q$ satisfies $|Q\cap\{w\ge t\}|\ge c_n\theta|Q|$.

For a dyadic subcube $Q'$ of $Q$ with center $z$ and side length $l$, let $s=l/4$. If $|Q'\cap\{w\ge\tau\}|>(1-\varepsilon')|Q'|$, where $\varepsilon'>0$ is chosen so that $\varepsilon' l^n\le\varepsilon|B(z,l/2)|$, {\color{black}the application of \eqref{eq:nondivpotentialfirstgrowth} with $R=s$ and $f=1-w/\tau$ gives $w\ge\tau/2$ on $B(z,s)$.} From this ball, a chain of balls of radius $s$ with fourfold dilates in $B(x,2R)$ reaches every point of the dyadic parent of $Q'$ (or of $Q$ if $Q'=Q$). The number of balls depends only on $n$. Thus $w\ge a\tau$ on that parent, with $a\in(0,1)$ depending only on $n$, $\nu$, $\|b\|_{L^n(D\cap V)}$, and $K_0$.

If $|Q\cap\{w\ge\tau\}|\le(1-\varepsilon')|Q|$, take the maximal dyadic subcubes of $Q$ with density greater than $1-\varepsilon'$. By differentiation, they cover $Q\cap\{w\ge\tau\}$ up to a null set. Their parents have density at most $1-\varepsilon'$. The maximal cubes among these parents are disjoint, cover the same set up to a null set, and lie in $\{w\ge a\tau\}$. Consequently,
\[
 |Q\cap\{w\ge\tau\}|\le(1-\varepsilon')|Q\cap\{w\ge a\tau\}|.
\]
Starting at $\tau=t$ and reducing the level successively by the factor $a$, after a number of steps depending only on $\theta$, $n$, and $\varepsilon'$, the density in $Q$ exceeds $1-\varepsilon'$. The estimate on $Q$ then gives a lower bound by a fixed positive multiple of $t$ on $Q$. Another chain of balls of radius $R/(400\sqrt n)$, with fourfold dilates in $B(x,2R)$, propagates the bound to $B(x,R/2)$; the number of balls depends only on $n$. Therefore, for every $\theta\in(0,1)$,
\begin{equation}\label{eq:nondivpotentialCD}
 |\{w\ge t\}\cap B(x,R)|\ge\theta|B(x,R)|
 \quad\Longrightarrow\quad
 \inf_{B(x,R/2)}w\ge\beta t,
\end{equation}
where $\beta>0$ depends only on $n$, $\nu$, $\|b\|_{L^n(D\cap V)}$, $K_0$, and $\theta$.
}For the balls in (CD),
\[
 R^2\|c\|_{L^{\infty}(B(x,2R))}
 \le\frac{\|\delta_D^2c\|_{L^{\infty}(D)}}{(A_{\mathrm{CD}}-2)^2}.
\]
Taking $\theta=\vartheta$ in~\eqref{eq:nondivpotentialCD}
proves (CD).

For (UC), let $g=v-tu$ and $T>0$. Since $L_cg\ge0$,
on the open set $\{g>T\}$ we have
\[
 L(g-T)\ge cg\ge0.
\]
{ The homogeneous first growth lemma, applied on these positive superlevel sets as in Section~\ref{sec:A2nondiv}, proves boundary (UC). Its threshold is independent of $c$.}
For a nonnegative strong solution $h$,
\[
 Lh=ch\ge0,
\]
{\color{black}so the homogeneous first growth lemma in Section~\ref{sec:A2nondiv}, applied to $h-\ell-T$ on $\{h>\ell+T\}$ for $\ell\ge0$ and $T>0$, gives interior (UC).}

Finally, put $M_c=\|\delta_D^2c\|_{L^{\infty}(D)}$.
{ On a testing neighborhood of radius $s$ where $s^2\|c\|_{L^\infty}\le\kappa_c$ (the norm is taken over that neighborhood), {\color{black}the critical-density estimate \eqref{eq:nondivpotentialCD} holds with $K_0=\kappa_c$} for every interior testing ball contained in the neighborhood, and (UC) is independent of $c$. Hence Proposition~\ref{prop:interior} supplies a Harnack constant $H_0\ge1$ and dilation $C_{\mathrm{int}}$, depending only on $n$, $\nu$, and $\|b\|_{L^n(D\cap V)}$.} For $B(x,4R)\subset D\cap V$,
\[
 \|c\|_{L^{\infty}(B(x,2R))}\le\frac{M_c}{4R^2}.
\]
Choose $\varepsilon>0$ so that
\[
 C_{\mathrm{int}}\varepsilon\le\frac12,
 \qquad
 \frac{C_{\mathrm{int}}^2\varepsilon^2}{4}\le\kappa_c,
\]
and set $\rho=\varepsilon R/(1+\sqrt{M_c})$.
For every $z\in B(x,R)$,
\[
 \begin{aligned}
 B(z,C_{\mathrm{int}}\rho)&\subset B(x,2R),\\
 (C_{\mathrm{int}}\rho)^2
 \|c\|_{L^{\infty}(B(z,C_{\mathrm{int}}\rho))}
 &\le \frac{C_{\mathrm{int}}^2\varepsilon^2M_c}
 {4(1+\sqrt{M_c})^2}\le\kappa_c.
 \end{aligned}
\]
Thus Proposition~\ref{prop:interior} applies with constant $H_0$
on all these balls. A chain of at most $C(1+\sqrt{M_c})$
balls joins any two points of $B(x,R)$, so
\[
 \sup_{B(x,R)}w
 \le H_0^{C(1+\sqrt{M_c})}\inf_{B(x,R)}w
 \le C_{\mathrm H}e^{\gamma_{\mathrm H}\sqrt{M_c}}
 \inf_{B(x,R)}w.
\]
Here $C_{\mathrm H}=H_0^C$ and
$\gamma_{\mathrm H}=C\max\{1,\log H_0\}$.
Theorems~\ref{thm:bhp} and~\ref{thm:quotient} give the
boundary comparison and quotient estimates.
\end{proof}

\endgroup

{ For fully nonlinear equations, boundary comparison in Lipschitz graph domains for constant-coefficient operators was established in \cite{Feldman01}; see also \cite{Ferrari06,ArgiolasFerrari09}.} { For Fatou theorems for fully nonlinear elliptic equations, see \cite[Section~II.3]{FGMS88}.}

\begin{example}[Pucci operators]\label{ex:pucci}
Let $X=\mathbb R^n$ with Lebesgue measure, $0<{\lambda_{\mathrm P}}\le{\Lambda_{\mathrm P}}$, and define, for the eigenvalues $e_j$ of a symmetric matrix $A$,
\[
 \mathcal M^+(A)={\Lambda_{\mathrm P}}\sum_{e_j>0}e_j+{\lambda_{\mathrm P}}\sum_{e_j<0}e_j,
 \qquad \mathcal M^-(A)={\lambda_{\mathrm P}}\sum_{e_j>0}e_j+{\Lambda_{\mathrm P}}\sum_{e_j<0}e_j.
\]
Take $F=\mathcal M^+$ or $F=\mathcal M^-$. More generally, let $F$ be continuous, positively homogeneous of degree one, independent of $x$, and satisfy
$\mathcal M^-(A-B)\le F(A)-F(B)\le\mathcal M^+(A-B)$.
Let $\Hcone$ be the continuous viscosity solutions of $F(D^2u)=0$, $\Scone$ its continuous viscosity supersolutions $F(D^2u)\le0$, and $\mathcal P$ the continuous viscosity supersolutions of $\mathcal M^-(D^2w)\le0$.
The viscosity difference principle gives $tu-v\in\mathcal P$ for $u\in\Scone$, $v\in\Hcone$, and $t\ge0$. Here positive homogeneity makes $tu$ a supersolution for $t>0$; the case $t=0$ follows from the extremal inequalities. Constants and positive rescalings preserve $\mathcal P$. Furthermore,
\[
 \mathcal C=\{{u}:\mathcal M^-(D^2{u})\le0\le\mathcal M^+(D^2{u})\}.
\]
The weak Harnack inequality for $\mathcal M^-$ gives (CD). For ${u}\ge0$ in $\mathcal C$, $({u}-c)_+$ is a viscosity subsolution of $\mathcal M^+\ge0$ for $c\ge0$. The local $L^2$ upper estimate gives interior (UC) by Section~\ref{sec:A2upper}. These standard viscosity estimates and the extremal difference principle are treated in \cite{CC95}.

For boundary comparison on an open neighborhood $V$, assume $u\ge0$ in $\Scone(D\cap V)$ and $v>0$ in $\Hcone(D\cap V)$, with $v$ extending continuously as zero on $\partial D\cap V$. For every $t\ge0$, ${w}=(tu-v)_-=(v-tu)_+$ is a viscosity subsolution of $\mathcal M^+\ge0$ inside $D\cap V$, and $0\le {w}\le v$. Its zero extension is continuous on $V$ and is a viscosity subsolution there. Indeed, at a boundary point an upper test function is nonnegative near its zero contact value, so its Hessian is positive semidefinite. The ambient local $L^2$ upper estimate therefore gives (UC), uniformly in $t$. { Its constant depends only on $n$, ${\lambda_{\mathrm P}}$, and ${\Lambda_{\mathrm P}}$, so neither the comparison coefficient $t$ nor the chosen pair enters the density constant.} Theorem~\ref{thm:bhp} applies on diameter-uniform domains. When both positive solutions vanish continuously on $\partial D\cap V$, (UC) also holds for $(tv-u)_-$, and Theorem~\ref{thm:quotient} applies.

{ For $n=2$, $\lambda_{\mathrm P}=1$, $\Lambda_{\mathrm P}=3$, and $F=\mathcal M^+$,} the positive functions $u=4+x^2-3y^2$ and $v=4+y^2-3x^2$ solve $F(D^2u)=F(D^2v)=0$ in the unit disk, but $F(D^2(u+v))=-8$.
\end{example}

\begin{corollary}[Pucci equations with critical drift and nonnegative potential]
\label{cor:pucci-critical-drift}
Let $X=\mathbb R^n$, $n\ge2$, with Lebesgue measure, and let
$D\subsetneq\mathbb R^n$ be diameter-uniform. Fix
$0<\lambda_{\mathrm P}\le\Lambda_{\mathrm P}$, and suppose that
$\beta\ge0$, $\beta\in L^n(D)$, $c\ge0$, and
$\delta_D^2c\in L^\infty(D)$. Fix either of the operators
\[
\begin{aligned}
 F_+[w]&=\mathcal M^+(D^2w)+\beta(x)|Dw|-c(x)w,\\
 F_-[w]&=\mathcal M^-(D^2w)-\beta(x)|Dw|-c(x)w.
\end{aligned}
\]
Let $\xi\in\partial D$ and $0<r<r_*/(4\Lambda_0)$, where
$\Lambda_0$ and $a_{\xi,r}$ are as in
Theorem~\ref{thm:bhp}. Put $V=B(\xi,\Lambda_0r)$.
Suppose $u,v\in W^{2,n}_{\mathrm{loc}}(D\cap V)\cap C(D\cap V)$
are positive strong solutions of the same equation $F_+[w]=0$
or $F_-[w]=0$ almost everywhere in $D\cap V$.
Assume that both extend continuously to $\overline D\cap V$
and vanish on $\partial D\cap V$. Then there exist
$C\ge1$ and $\alpha\in(0,1)$ such that, for every $x\in D\cap B(\xi,r)$,
\begin{equation}\label{eq:pucci-critical-bhp}
 C^{-1}\frac{u(a_{\xi,r})}{v(a_{\xi,r})}
 \le\frac{u(x)}{v(x)}\le
 C\frac{u(a_{\xi,r})}{v(a_{\xi,r})}.
\end{equation}
Moreover, for every $x,y\in D\cap B(\xi,r/16)$,
\begin{equation}\label{eq:pucci-critical-holder}
 \left|\frac{u(x)}{v(x)}-\frac{u(y)}{v(y)}\right|
 \le C\left(\frac{|x-y|}{r}\right)^\alpha
 \frac{u(a_{\xi,r})}{v(a_{\xi,r})}.
\end{equation}
The constants depend only on $n$, $\lambda_{\mathrm P}$,
$\Lambda_{\mathrm P}$, the diameter-uniform constant of $D$,
$\|\beta\|_{L^n(D)}$, and
$\|\delta_D^2c\|_{L^{\infty}(D)}$.
\end{corollary}

\begin{proof}
Fix one of the two operators $F=F_+$ or $F=F_-$.
On every open $U\subset D$, take $\Scone(U)$ to be the
continuous $W^{2,n}_{\mathrm{loc}}$ strong supersolutions
of $F[w]\le0$, and $\Hcone(U)$ the corresponding strong
solutions. Define the comparison supersolution family by
\[
 \mathcal P(U)=\bigl\{z\in W^{2,n}_{\mathrm{loc}}(U)\cap C(U):
 \mathcal M^-(D^2z)-\beta|Dz|-cz\le0\bigr\}.
\]
As in Example~\ref{ex:pucci}, the extremal difference inequality
and positive homogeneity give, almost everywhere in $U$,
\[
 \mathcal M^-(D^2(tw-h))-\beta|D(tw-h)|-c(tw-h)
 \le tF[w]-F[h]\le0
\]
for $w\in\Scone(U)$, $h\in\Hcone(U)$, and $t\ge0$.
{ Since $c\ge0$,}
$\mathcal P(U)$ is closed under nonnegative constant shifts
and positive scalar multiplication. Thus (Comp) holds.

For a symmetric matrix $H$, write $\Pi_-(H)$ and $\Pi_+(H)$
for the orthogonal spectral projections onto its negative and
positive eigenspaces. Set
\[
 \begin{aligned}
 A^-(H)&=\lambda_{\mathrm P}I
 +(\Lambda_{\mathrm P}-\lambda_{\mathrm P})\Pi_-(H),\\
 A^+(H)&=\lambda_{\mathrm P}I
 +(\Lambda_{\mathrm P}-\lambda_{\mathrm P})\Pi_+(H).
 \end{aligned}
\]
These are Borel functions of $H$ and satisfy
\[
 \operatorname{tr}(A^\pm(H)H)=\mathcal M^\pm(H),
 \qquad
 \lambda_{\mathrm P}I\le A^\pm(H)\le\Lambda_{\mathrm P}I.
\]
For every nonnegative $z\in\mathcal P(U)$, take
$A_z(x)=A^-(D^2z(x))$ and
\[
 q_z(x)=
 \begin{cases}
 -\beta(x)Dz(x)/|Dz(x)|,&Dz(x)\ne0,\\
 0,&Dz(x)=0.
 \end{cases}
\]
Then, almost everywhere in $U$,
\[
 \operatorname{tr}(A_zD^2z)+q_z\cdot Dz-cz\le0,
 \qquad |q_z|\le\beta.
\]
The ellipticity hypothesis of Theorem~\ref{thm:nonnegativepotentials}
holds with $\nu=\min\{1,\lambda_{\mathrm P},\Lambda_{\mathrm P}^{-1}\}$.
If $B(x,A_{\mathrm{CD}}s)\subset U$, then
\[
 \|q_z\|_{L^n(B(x,2s))}\le\|\beta\|_{L^n(D)},\qquad
 s^2\|c\|_{L^{\infty}(B(x,2s))}
 \le\frac{\|\delta_D^2c\|_{L^{\infty}(D)}}
 {(A_{\mathrm{CD}}-2)^2}.
\]
The linear critical-density estimate
\eqref{eq:nondivpotentialCD}, proved in
Theorem~\ref{thm:nonnegativepotentials}, therefore applies
uniformly to these measurable coefficients and proves (CD).

{ If $h\ge0$ belongs to $\mathcal C(U)$, then $\mathcal M^+(D^2h)+\beta|Dh|-ch\ge0$. Linearize the first two terms using $A^+(D^2h)$ and the vector $\beta Dh/|Dh|$ on $\{Dh\ne0\}$, with vector zero elsewhere. The resulting linear operator satisfies $L_h(h-\ell-T)\ge ch\ge0$ on $\{h>\ell+T\}$ for $\ell\ge0$ and $T>0$. {\color{black}The ellipticity and drift bounds are $\nu=\min\{1,\lambda_{\mathrm P},\Lambda_{\mathrm P}^{-1}\}$ and $\|\beta\|_{L^n(D)}$, respectively.} Thus the first growth lemma in Section~\ref{sec:A2nondiv} gives interior (UC).

For the boundary estimate, take a nonnegative strong supersolution $u$ and a positive strong solution $v$ extending continuously as zero on $\partial D\cap V$. The comparison-difference inequality gives $\mathcal M^+(D^2g)+\beta|Dg|-cg\ge0$ for $g=v-tu$ and $t\ge0$. The same linearization, now for $g$, yields $L_g(g-T)\ge cg\ge0$ on $\{g>T\}$ for $T>0$. Since $g_+\le v$, Section~\ref{sec:A2nondiv} applies on these positive superlevel domains without requiring $g_+\in W^{2,n}_{\mathrm{loc}}$. This gives boundary (UC) uniformly in $t$, and the first-growth threshold is independent of $c$. Interchanging two positive solutions with the stated boundary condition gives the reverse estimate. }Theorems~\ref{thm:bhp} and
\ref{thm:quotient} yield \eqref{eq:pucci-critical-bhp}
and \eqref{eq:pucci-critical-holder}.
\end{proof}

\section{Examples illustrating the limits of the axioms}\label{sec:limits}

\begin{example}[Affine functions without minimum closure]\label{ex:sparse}
\begingroup
{ On the space and domain of Example~\ref{ex:concave}, let $\Scone(U)=\Hcone(U)=\mathcal P(U)$ consist of functions affine on each component of $U$. Then (Comp) and (CD) hold, and the interior and boundary (UC) estimates under the conditions in that example follow by convexity.}
{\color{black}The functions $x$ and $1-x$ belong to $\Scone(D)$, whereas their minimum does not.}
\endgroup
\end{example}

\begin{example}[{Constant functions with nonzero boundary values}]\label{ex:vacuous}
\begingroup
{ With the same $X$ and $D$, let $\Scone(U)=\Hcone(U)=\mathcal P(U)$ consist only of constant functions. Then (Comp), (CD), and interior (UC) hold.}
{ Near either endpoint, (UC) holds on balls centered in $D$ with $M>2$. Indeed, if the constant value of the nonnegative comparison function exceeds the testing level, the superlevel set occupies at least half of the ambient interval whenever it meets only that endpoint. Thus the density hypothesis in \eqref{eq:UD} implies that the function is at most the testing level. Positive solutions in this example are nonzero constants and need not vanish at the boundary.}
\endgroup
\end{example}

\begingroup

\section{Verification of (UC)}\label{sec:verifyA2}
{Standard local upper estimates imply (UC) for the functions ${w}=({tu}-v)_-$ required in boundary comparison, uniformly for ${t}\ge0$, including ${w}=v$ when ${t}=0$. The reverse comparison requires the analogous estimates with $u$ and $v$ interchanged.}

\subsection{The local upper estimate}\label{sec:A2upper}
Let ${f}:D\cap V\to[0,\infty)$ be measurable. Assume that, for every compact set $K\subset V$,
\[
 \int_{D\cap K}{f}^2\,d\mu<\infty.
\]
Suppose that
\begin{equation}\label{eq:A2localupper}
 \esssup_{D\cap B(x,s)}{f}
 \le C_0\left(\frac{1}{\mu(B(x,2s))}
                   \int_{D\cap B(x,2s)}{f}^2\,d\mu\right)^{1/2}{,}
\end{equation}
whenever $\overline{B(x,2s)}\subset V$, with a common $C_0\ge1$. A finite covering by smaller balls gives the boundedness required in Definition~\ref{def:UC}. For $B(x,4s)\subset V$, the local upper estimate \eqref{eq:A2localupper} and the finite-cover argument give finite essential supremum on $D\cap B(x,2s)$. Under the density hypothesis in \eqref{eq:UD}, splitting at level $t$ yields
\begin{gather*}
 \frac{1}{\mu(B(x,2s))}\int_{D\cap B(x,2s)}{f}^2\,d\mu
 \le t^2+M^{-1}{\bigl(\esssup_{D\cap B(x,2s)}{f}\bigr)^2},\\
 \esssup_{D\cap B(x,s)}{f}\le C_0t+C_0M^{-1/2}{\esssup_{D\cap B(x,2s)}{f}}.
\end{gather*}
{Thus (UC) follows after increasing the constant $M$ in (CD), if necessary, so that $M\ge\max\{2,C_0,4C_0^2\}$. Increasing $M$ preserves (CD).} The denominator in \eqref{eq:A2localupper} is the measure of the entire ambient ball.

\subsection{Divergence-form equations and Dirichlet forms}\label{sec:A2div}
Consider the operator in Example~\ref{ex:elliptic}, with the same
coefficient assumptions, locally on $D\cap V$, where
$D\subsetneq\mathbb R^n$ and $V\subset\mathbb R^n$ are open.
Write $\lambda_{\mathrm{ell}}$ and $\Lambda_{\mathrm{ell}}$ for its
ellipticity and boundedness constants. Its interior local upper estimate
\eqref{eq:A2localupper} gives (UC) by Section~\ref{sec:A2upper}.

Let $u\ge0$ be a weak supersolution and $v\ge0$ a weak solution in
$D\cap V$. For ${t}\ge0$, the function ${w}=(v-{tu})_+$ is a nonnegative
weak subsolution with ${w}\le v$. For boundary balls it suffices to show, for every $\eta\in C_c^\infty(V)$,
\begin{equation}\label{eq:A2zerotrace}
 \eta {w}\in H_0^1(D),
\end{equation}
with localized functions extended by zero outside $V$.
Indeed, testing with $\phi T_\varepsilon({w})$, where
$0\le\phi\in C_c^\infty(V)$ and
$T_\varepsilon(t)=\min\{t/\varepsilon,1\}$, gives
\[
 \int_{D\cap V}T_\varepsilon({w})\bm{A}\nabla {w}\cdot\nabla\phi\,dx
 +\frac1\varepsilon
 \int_{D\cap V\cap\{0<{w}<\varepsilon\}}
 \phi\bm{A}\nabla {w}\cdot\nabla {w}\,dx\le0.
\]
The second term is nonnegative. Letting $\varepsilon\downarrow0$,
using $\nabla {w}=0$ almost everywhere on $\{{w}=0\}$, shows that the
zero extension of ${w}$ is a weak subsolution in $V$, after extending
$\bm{A}$ elliptically outside $D$. Thus \eqref{eq:A2localupper}
gives (UC) on $V$, uniformly in ${t}$.

\begin{lemma}[Finite energy of a positive part]\label{lem:A2traceenergy}
Let $u\ge0$ be a weak supersolution in $D\cap V$ and let
$v\ge0$ belong to $H^1_{\mathrm{loc}}(D\cap V)$. If
$\eta v\in H_0^1(D)$ for every nonnegative
$\eta\in C_c^\infty(V)$, then
$\eta(v-u)_+\in H_0^1(D)$ for every $\eta\in C_c^\infty(V)$.
\end{lemma}
\begin{proof}
Fix $\eta$ and choose $0\le\chi\in C_c^\infty(V)$ equal to one
near $\operatorname{supp}\eta$. Take nonnegative
$f_j\in C_c^\infty(D\cap V)$ converging to $\chi v$ in $H^1(D)$,
and set ${w}_j=(f_j-u)_+$. Then ${w}_j\in H_0^1(D)$ and
$0\le {w}_j\le f_j$. Testing the supersolution inequality with ${w}_j$
gives
\[
\|\nabla {w}_j\|_{L^2(D)}
 \le \frac{ \Lambda_{\mathrm{ell}} }{ \lambda_{\mathrm{ell}} }  \|\nabla f_j\|_{L^2(D)}.
\]
Thus $({w}_j)$ is bounded in $H_0^1(D)$ and converges in $L^2(D)$
to $(\chi v-u)_+$. Weak closedness and multiplication by $\eta$
give the assertion.
\end{proof}

{For $t,c\ge0$, the function $tu+c$ is a nonnegative weak supersolution. Thus Lemma~\ref{lem:A2traceenergy} gives $\eta(v-tu-c)_+\in H_0^1(D)$ for every $\eta\in C_c^\infty(V)$.} The argument also applies to the scalar
$A_2$-weighted form in Example~\ref{ex:spectral}, with its weighted
zero-trace space and local upper estimate.

Alternatively, \eqref{eq:A2zerotrace} follows if, for every compact set $K\subset V$,
\begin{equation}\label{eq:A2HV}
 \int_{D\cap K}\left(v^2+\frac{v^2}{\del^2}\right)\,dx<\infty.
\end{equation}
Indeed, let $\theta_\varepsilon$ vanish on $\{\del\le\varepsilon\}$
and equal one on $\{\del\ge2\varepsilon\}$, with
$|\nabla\theta_\varepsilon|\le C/\varepsilon$.
The Caccioppoli inequality for ${w}$ gives
\[
 \int_{D\cap V}\eta^2\theta_\varepsilon^2|\nabla {w}|^2\,dx
 \le C\int_{D\cap V}{w}^2
 (\theta_\varepsilon^2|\nabla\eta|^2
  +\eta^2|\nabla\theta_\varepsilon|^2)\,dx.
\]
By ${w}\le v$ and \eqref{eq:A2HV}, the right-hand side is uniformly
bounded; its boundary-cutoff term tends to zero. Hence
$\eta\theta_\varepsilon {w}$ is bounded in $H_0^1(D)$ and converges
to $\eta {w}$ in $L^2(D)$. Weak closedness proves
\eqref{eq:A2zerotrace}. The constant in (UC) depends only on the
structural coefficients, not on the integral in \eqref{eq:A2HV}.

\subsection{Variational Neumann boundary conditions}\label{sec:A2neumann}

We work with a variational formulation because the boundary of a
Euclidean uniform domain need not admit a normal vector. Let
$D\subsetneq\mathbb R^n$, $n\ge2$, be a Euclidean uniform domain with
uniformity constant $A_D$, and let $V\subset\mathbb R^n$ be open. If
$D$ is bounded, assume
\begin{equation}\label{eq:A2neumannscale}
 \operatorname{diam}V\le\frac14\operatorname{diam}D.
\end{equation}
We use the divergence-form operator of Example~\ref{ex:elliptic} and the ellipticity and boundedness constants $\lambda_{\mathrm{ell}}$ and $\Lambda_{\mathrm{ell}}$ fixed in Section~\ref{sec:A2div}.
We say that $h$ has locally finite
energy up to the boundary in $V$ if $\eta h\in H^1(D)$ for every
$\eta\in C_c^\infty(V)$. A variational Neumann weak supersolution of
$Lh\ge0$ has this property and satisfies
\begin{equation}\label{eq:A2neumannweak}
 \int_{D\cap V}\bm{A}\nabla h\cdot\nabla\varphi\,dx\ge0
\end{equation}
for every nonnegative $\varphi\in H^1(D)$ supported in $D\cap K$ for
some compact $K\Subset V$. Equality for all such test functions of
either sign defines a homogeneous variational Neumann weak solution.
In particular, the test functions need not vanish near $\partial D\cap V$.

\begin{proposition}[Neumann boundary verification of (UC)]\label{prop:A2neumann}
Under the preceding assumptions, let $u\ge0$ be a variational Neumann
weak supersolution and let $v\ge0$ be a homogeneous variational Neumann
weak solution in $D\cap V$. For $t,c\ge0$ set $w_{t,c}=(v-tu-c)_+$.
There exists $C_0\ge1$, depending only on $n$, $A_D$,
$\lambda_{\mathrm{ell}}$, and $\Lambda_{\mathrm{ell}}$, such that
\begin{equation}\label{eq:A2neumannupper}
 \esssup_{D\cap B(x,s)}w_{t,c}
 \le C_0\left(\frac1{|B(x,2s)|}
       \int_{D\cap B(x,2s)}w_{t,c}^2\,dy\right)^{1/2}
\end{equation}
for every $x\in D$ and $s>0$ with $\overline{B(x,2s)}\subset V$.
Consequently, the family $\{w_{t,c}:t,c\ge0\}$ satisfies (UC) on $V$
with a common constant independent of $t$, $c$, $u$, and $v$.
If both $u$ and $v$ are {\color{black}nonnegative} homogeneous variational Neumann weak
solutions, the same conclusion holds for $(u-tv-c)_+$.
\end{proposition}

\begin{proof}
The boundary variational inequalities imply that $v-tu-c$ is a weak
subsolution with boundary-touching test functions. Testing with a
nonnegative test function multiplied by a Lipschitz approximation to
$\mathbf1_{\{v-tu-c>0\}}$ and passing to the limit shows that
$w_{t,c}$ is also a weak subsolution in this sense.

For $\overline{B(x,2s)}\subset V$, the scale restriction above gives
$s\le\operatorname{diam}D/16$ if $D$ is bounded. Rescaling the uniform
domain and applying the Sobolev extension theorem
\cite[Theorem~1]{Jones81} to localized cutoffs yields a Sobolev
inequality with a constant independent of $x$ and $s$. The standard
boundary Caccioppoli inequality and Moser iteration, using the same
boundary-touching tests, then give \eqref{eq:A2neumannupper}.
Section~\ref{sec:A2upper} yields (UC) after enlarging the common
constant in (CD) so that $M\ge\max\{2,C_0,4C_0^2\}$. The same
argument applies with $u$ and $v$ interchanged when both are
Neumann weak solutions.
\end{proof}

{
Neumann heat kernel estimates on (inner) uniform domains have been
established in \cite{GS11,Mur24,Ant25}. {\color{black}As in \cite{Mur24}, the proof of Proposition~\ref{prop:A2neumann} uses the extension property, but does not require symmetry of the coefficient matrix.}}

\begin{remark}
The abstract BHP also applies to homogeneous Robin boundary conditions with nonnegative coefficients, under suitable assumptions on boundary traces.
{
For one-sided NTA pairs of mixed dimension, \cite[Theorem~6.6]{DDEFMM25}
gives boundary comparison for two nonnegative solutions of the same
homogeneous Robin problem with a constant Robin coefficient.
The comparison constant is independent of the Robin coefficient and the
scale. }
\end{remark}

\subsection{Hardy inequalities and a portion of the boundary}\label{sec:A2partialHardy}
{
We characterize the local zero-trace condition in Section~\ref{sec:A2div}
by a capacity weight on bounded Sobolev extension domains.}

For partial boundary data, let $F\subset\partial\Omega$ be closed and set
\[
 H_F^1(\Omega)
 =\overline{\{\varphi|_\Omega:\varphi\in C_c^\infty(\mathbb R^n),
                    \ \operatorname{supp}\varphi\cap F=\varnothing\}}^{\,H^1(\Omega)}.
\]

\begingroup
For a closed set $F\subset\partial\Omega$ and $\rho>0$, let
$\operatorname{cap}_2(K,G)$ denote the variational relative $2$-capacity
of $K\subset G$ and define
\begin{equation}\label{eq:A2capweight}
 \begin{split}
 c_{8,n}&=\operatorname{cap}_2(B(0,1),B(0,8)),\\
 \Theta_F(x,\rho)
 &=\frac{\operatorname{cap}_2(F\cap B(x,\rho),B(x,8\rho))}
              {c_{8,n}\rho^{n-2}},\\
 \mathcal V_F(x)&=2\int_0^\infty\Theta_F(x,\rho)^2
                  \frac{d\rho}{\rho^3}.
 \end{split}
\end{equation}
Here $c_{8,n}>0$ also for $n=2$, when the relative capacity is
scale invariant. We use outer relative capacity for noncompact sets.
If $F=\varnothing$, set $\mathcal V_F=0$. Otherwise,
$0\le\mathcal V_F(x)\le\operatorname{dist}(x,F)^{-2}$.

\begin{proposition}[Capacity-integrated characterization of partial zero trace]
\label{prop:A2capacitary}
Suppose $n\ge2$, $\Omega\subset\mathbb R^n$ is a bounded
$H^1$-extension domain, and $F\subset\partial\Omega$ is closed.
Then
\begin{equation}\label{eq:A2captrace}
 H_F^1(\Omega)=\left\{f\in H^1(\Omega):
     \int_\Omega\mathcal V_F(x)|f(x)|^2\,dx<\infty\right\}.
\end{equation}
Moreover, for every $f\in H_F^1(\Omega)$,
\begin{equation}\label{eq:A2capHardy}
 \int_\Omega\mathcal V_F(x)|f(x)|^2\,dx
 \le C\|f\|_{H^1(\Omega)}^2,
\end{equation}
where $C$ depends on $n$, $\Omega$, and its extension constant,
but not on $F$.
\end{proposition}
\begin{proof}
See Appendix~\ref{app:capacity}.
\end{proof}

\begin{corollary}[A capacitary boundary criterion for (UC)]
\label{cor:A2capacitary}
Suppose $D\subset\mathbb R^n$ is a bounded Euclidean uniform domain,
$n\ge2$, and $V\subset\mathbb R^n$ is an open boundary neighborhood.
Set $F=\partial D\cap\overline V$. Consider the elliptic operator
and comparison functions $u$ and $v$ in
Section~\ref{sec:A2div}, where $u\ge0$ is a weak supersolution and
$v\ge0$ is a weak solution in $D\cap V$.
{
The condition $\eta v\in H_0^1(D)$ for every $\eta\in C_c^\infty(V)$
is equivalent to the following condition for every nonnegative
$\eta\in C_c^\infty(V)$:
\begin{equation}\label{eq:A2capcondition}
 \eta v\in H^1(D),\qquad
 \int_D\mathcal V_F(x)|\eta(x)v(x)|^2\,dx<\infty.
\end{equation}
Under either condition, $\eta(v-tu)_+\in H_0^1(D)$ for every
$\eta\in C_c^\infty(V)$ and $t\ge0$, and (UC) holds on $V$ for
$(v-tu)_+$ with constants independent of $t$ and $F$.}

\end{corollary}
\begin{proof}
{
Uniformity gives the $H^1$-extension property of $D$. For a nonnegative
cutoff $\eta$, Proposition~\ref{prop:A2capacitary} identifies
\eqref{eq:A2capcondition} with $\eta v\in H_F^1(D)$.
Choose $\chi\in C_c^\infty(V)$ equal to one near $\operatorname{supp}\eta$.
If $\varphi_j|_D\to\eta v$ in $H^1(D)$ with
$\varphi_j\in C_c^\infty(\mathbb R^n\setminus F)$, then
$(\chi\varphi_j)|_D\in C_c^\infty(D)$ and these functions converge
to $\eta v$ in $H^1(D)$. Thus $\eta v\in H_0^1(D)$.
For a signed cutoff, multiply by a nonnegative cutoff equal to one
near its support. The converse follows from $H_0^1(D)\subset H_F^1(D)$.
Lemma~\ref{lem:A2traceenergy} applied to $tu$ gives the assertion for
$(v-tu)_+$, and Section~\ref{sec:A2div} gives (UC).}

\end{proof}
\endgroup

{
{\color{black}Hardy inequalities with boundary-distance weights} help verify
\eqref{eq:A2HV}. If $\mathbb R^n\setminus\Omega$ is uniformly $2$-fat,
then \cite[Corollary~3.11]{KM97} yields, for every $f\in H_0^1(\Omega)$,
\begin{equation}\label{eq:A2KMHardy}
 \int_\Omega\frac{|f|^2}{\operatorname{dist}(x,\partial\Omega)^2}\,dx
 \le C\int_\Omega|\nabla f|^2\,dx.
\end{equation}
Conversely, for any open $\Omega$ and $f\in H^1(\Omega)$,
\cite[Theorem~3.13]{KM97} gives
\begin{equation}\label{eq:A2KMtrace}
 \frac{f}{\operatorname{dist}(\cdot,\partial\Omega)}\in L^2(\Omega)
 \quad\Longrightarrow\quad f\in H_0^1(\Omega).
\end{equation}
Section~\ref{sec:A2div} supplies the requisite energy bound for the positive part.

For closed $F\subset\partial\Omega$ and $R>0$, suppose that every
$f\in H_F^1(\Omega)$ satisfies
\begin{equation}\label{eq:A2partialHardy}
 \int_\Omega\frac{|f|^2}{d(x,F)^2}\,dx
 \le C_F\left(\int_\Omega|\nabla f|^2\,dx
                +R^{-2}\int_\Omega|f|^2\,dx\right).
\end{equation}
For ambient Lipschitz functions vanishing on uniformly $p$-fat sets,
see \cite[Theorem~1.3 and Proposition~5.4]{LTV17} and
\cite[Chapters~6, 7, and~10]{KLV21}; for thickness, Poincar\'e, and
trace-preserving extension conditions on $\Omega$, see
\cite[Theorem~3.1]{EHR15}. For barrier and two-weight criteria, see
\cite[Corollary~5.3]{Har24} and \cite[Theorem~3.11 and Corollary~3.12]{Har26}.

Fix $\xi\in\partial D$ and $0<r<R'<R$ with $\overline{B(\xi,R)}\subset V$.
Set $\Omega=D\cap B(\xi,R)$ and let $F\subset\partial\Omega\cap\partial D$
be closed and contain $\partial D\cap\overline{B(\xi,R')}$.
On $D\cap B(\xi,r)$, $d(\cdot,F)=\del$ whenever $\del<R'-r$;
thus $\del^{-1}\le\max\{d(\cdot,F)^{-1},(R'-r)^{-1}\}$.
Hence, if $\eta v\in H_F^1(\Omega)$ for every $\eta\in C_c^\infty(B(\xi,r))$,
\eqref{eq:A2partialHardy} yields \eqref{eq:A2HV} locally.}

\begin{remark}[Potential-theoretic zero boundary values]\label{rem:A2qe}
{
For bounded harmonic $v$ in $D\cap V$, the condition
$\eta v\in H_0^1(D)$ for all $\eta\in C_c^\infty(V)$ is equivalent to
$v$ tending to zero at every Wiener regular point of $\partial D\cap V$.
This local equivalence requires no global $H^1(D)$ assumption.}

{
For example, let $D=B(0,1)\setminus\{0\}\subset\mathbb R^2$,
$V=B(0,1/2)$, and $v=1$. The point $0$ is polar, so
$\mathcal V_{\{0\}}=0$ and $\eta v\in H_0^1(D)$ for every
$\eta\in C_c^\infty(V)$, while
$\int_{D\cap B(0,1/4)}v^2/\del^2\,dx=\infty$.

For the Laplacian, a nonnegative subharmonic function locally bounded up
to a boundary patch with upper limit zero off a polar subset of the patch
has a subharmonic extension across it by zero extension and upper
regularization. See \cite[Section~5]{AHL06} for the harmonic case.
For $(v-tu)_+$, domination by $v$ gives the required boundedness and
boundary limits, and the ambient mean-value inequality gives
\eqref{eq:A2localupper}.}

\end{remark}

\subsection{Boundary-continuous solutions of nondivergence equations}\label{sec:A2nondiv}
For the operator and functions of Example~\ref{ex:nondivlinear}, the boundary estimate follows from \cite[Lemma~2.1]{Safonov} on positive superlevel domains. The scale-invariant drift estimate in Example~\ref{ex:nondivlinear} applies on each rescaled testing neighborhood.

For ${w}=(v-tu)_+$, $t\ge0$, a level ${T}>0$, and a ball $B(x,4s)\subset V$, apply the cited lemma to $f=v-tu-{T}$ on
\[
 \Omega_{{T}}=\{v-{tu}>{T}\}\cap D\cap B(x,3s).
\]
Since $u\ge0$ and $t\ge0$, every point of $\Omega_{{T}}$ satisfies $v>T$. Moreover, $\overline{B(x,3s)}\subset V$ and $v$ extends continuously to $\overline D\cap V$, with $v=0$ on $\partial D\cap V$. Consequently,
\[
 \overline{\Omega_{{T}}}\cap\partial D\cap\overline{B(x,3s)}=\varnothing.
\]
Here $Lf\ge0$, and $f$ vanishes on the level-set boundary inside $D$. Therefore the boundary hypothesis of \cite[Lemma~2.1]{Safonov} holds without any boundary-continuity assumption on $u$. With its contraction factor fixed at $1/2$, that lemma gives a threshold $\varepsilon_0>0$, depending only on $n$, $\nu$, and $\|b\|_{L^n(D\cap V)}$, such that
\[
 \frac{|D\cap B(x,2s)\cap\{{w}>{T}\}|}{|B(x,2s)|}\le\varepsilon_0
 \quad\Longrightarrow\quad
 \sup_{D\cap B(x,s)}{w}\le{T}+{\frac{1}{2}}\sup_{D\cap B(x,2s)}{w}.
\]
{This is (UC) after increasing the constant $M$ in (CD), if necessary, so that $M\ge\max\{2,\varepsilon_0^{-1}\}$.} Local boundedness follows from ${w}\le v$ and the continuous zero boundary values. For a nonnegative strong solution ${u}\in W^{2,n}_{\mathrm{loc}}(U)\cap C(U)$ with $L{u}=0$ in $U\subset D\cap V$, the interior truncation $({u}-c)_+$, $c\ge0$, is treated on $B(x,4s)\subset U$ by taking $f={u}-c-{T}$, $T>0$. Only interior level-set boundaries occur. Interchanging two positive solutions with the same zero boundary values gives (UC) for $(tv-u)_-$.

\begingroup
\begin{remark}[Boundary regularity for nondivergence equations]\label{rem:nondivboundaryregularity}
For boundary regularity of linear nondivergence equations, see \cite{Kry66,Mi70}; for boundary comparison with continuous zero boundary values, see \cite{Bauman84}. {\color{black}For boundary regularity of the Pucci extremal operators, see the sufficiency direction of the Wiener criterion in \cite[Theorem~1.2]{LeeLee23}.}
\end{remark}
\par\endgroup
\par\endgroup

\section*{Acknowledgments}
{ The author thanks Professor Hiroaki Aikawa for his comments at a research meeting in February 2026, where the author presented an argument deriving almost positivity from porosity estimates.}
This work was supported by JSPS KAKENHI Grant Numbers 21K03281 and 23H03798.
The author used ChatGPT (OpenAI) to assist with mathematical discussions,
manuscript preparation, and English-language editing.

\appendix
\begingroup
\section{Proof of the capacity-integrated zero-trace criterion}\label{app:capacity}

Chen, Xiao, Yang, Yuan, and Zhang \cite{ChenXYYZ26} solved
Maz'ya's open Problem~8 on a capacitary-distance Hardy inequality
using probabilistic ideas based on killed Brownian motion. The
capacity-integrated Hardy inequality proved in this appendix was
inspired by their result.

This appendix proves Proposition~\ref{prop:A2capacitary} for every $n\ge2$.
In dimension $n=2$, relative capacity is used throughout.

\subsection{Relative capacity and a full Dirichlet Hardy inequality}

We use the relative capacity $\operatorname{cap}_2$ and the weights $\Theta_E$ and $\mathcal V_E$ from \eqref{eq:A2capweight}, replacing $F$ by $E$. The capacity is interpreted through its equivalent Sobolev obstacle formulation with quasi-everywhere constraints; outer regularization applies to arbitrary sets.
These quantities are measurable by outer regularity and compact
approximation. Moreover, $\Theta_E(x,R)\le1$, and
$\Theta_E(x,R)=0$ for $R<\operatorname{dist}(x,E)$.

\begin{lemma}[Comparison of concentric condenser capacities]
\label{lem:appcondenser}
Fix $1<a<b<c<\infty$. For $K\subset B(x,aR)$ and a fixed
$\varepsilon>0$ with $a+\varepsilon<b$, one has
\begin{equation}\label{eq:appcondenser}
 \operatorname{cap}_2(K,B(x,bR))
 \le C\operatorname{cap}_2(K,B(x,cR)),
\end{equation}
{\color{black}where $C$ depends only on $n$, $a$, $b$, $c$, and $\varepsilon$.}
In particular, if $|x-y|\le R$ and $K\subset B(y,R)$, then
$\operatorname{cap}_2(K,B(y,8R))\le C_n\operatorname{cap}_2(K,B(x,8R))$.
\end{lemma}
\begin{proof}
Choose a smooth cutoff equal to one on $B(x,(a+\varepsilon/2)R)$
and supported in $B(x,bR)$. Multiply an admissible function in
$H_0^1(B(x,cR))$ by this cutoff. The Poincar\'e inequality on the
larger ball controls the cutoff term by the original Dirichlet energy,
uniformly under scaling. For the shifted-center estimate, multiply an admissible
function on $B(x,8R)$ by a cutoff supported in $B(y,4R)$ and equal to
one on $B(y,2R)$; the Poincar\'e inequality on $B(x,8R)$ controls the
cutoff term. Outer regularization gives the assertion for general sets.
\end{proof}

\begin{lemma}[Full Dirichlet estimate]\label{lem:appfullhardy}
Let $G\subset\mathbb R^n$ be any bounded open set, $n\ge2$, and
$E=\mathbb R^n\setminus G$. Then, for every $g\in H_0^1(G)$,
\begin{equation}\label{eq:appfullhardy}
 \int_G\mathcal V_E(x)|g(x)|^2\,dx
 \le C_n\int_G|\nabla g|^2\,dx.
\end{equation}
\end{lemma}
\begin{proof}
Let $A_G$ be the nonnegative Dirichlet Laplacian on $L^2(G)$,
associated with the form $\mathcal E(g)=\int_G|\nabla g|^2\,dx$.
For $R>0$ and $f\in L^2(G)$, put
\[
 J_R f=(I+R^2 A_G)^{-1}f.
\]
This is the unique minimizer of
$\mathcal E(v)+R^{-2}\|f-v\|_{L^2(G)}^2$ in $H_0^1(G)$.
By order preservation, $0\le J_Rf\le J_R1\le1$ whenever
$0\le f\le1$. For every $g\in H_0^1(G)$, the spectral theorem gives the exact identity
\begin{equation}\label{eq:appscale}
 2\int_0^\infty\|g-J_Rg\|_{L^2(G)}^2\frac{dR}{R^3}
 =\mathcal E(g).
\end{equation}
Indeed, on each spectral value $\lambda\ge0$, the multiplier on
the left integrates to
$2\int_0^\infty (R^2\lambda/(1+R^2\lambda))^2
\,dR/R^3=\lambda$.

Set $z_R=1-J_R1$. It takes values in $[0,1]$, equals one
quasi-everywhere on $E$ after extending $J_R1$ by zero, and satisfies, in $G$,
\begin{equation}\label{eq:appabsorption}
 -\Delta z_R+R^{-2}z_R=0.
\end{equation}
We claim that, for every $R>0$ and $x\in G$, the following pointwise capacitary lower bound holds:
\begin{equation}\label{eq:appobstacle}
 z_R(x)\ge c_n\Theta_E(x,R).
\end{equation}
It suffices first to work with a compact set
$K\subset E\cap B(x,R)$. Let $h$ minimize
\[
 \int_{B(x,8R)}\bigl(|\nabla h|^2+R^{-2}h^2\bigr)\,dy
\]
over $H_0^1(B(x,8R))$ with $h\ge1$ quasi-everywhere on $K$.
Truncation yields $0\le h\le1$, with $h=1$ quasi-everywhere on
$K$. It is a supersolution of the equation in
\eqref{eq:appabsorption} and a solution away from $K$.
Write $M=h(x)$; the value is well-defined using interior
continuity since $x\notin K$. After rescaling to $R=1$, the
interior weak Harnack inequality for nonnegative supersolutions,
with exponent one, gives
\[
 \int_{B(x,4R)}h\,dy\le C_n MR^n.
\]
The exponent one is admissible for the uniformly elliptic
absorption equation in every dimension $n\ge2$.
Interior local boundedness and the Caccioppoli inequality in the
annulus $B(x,7R/2)\setminus B(x,3R/2)$, which avoids $K$,
then give
\[
 \|h\|_{L^{\infty}(B(x,3R)\setminus B(x,2R))}
 \le C_n M,\qquad
 \int_{B(x,3R)\setminus B(x,2R)}|\nabla h|^2\,dy
 \le C_n M^2 R^{n-2}.
\]
Choose $\eta\in C_c^\infty(B(x,3R))$ equal to one on
$B(x,2R)$ with $|\nabla\eta|\le C/R$. The test
$\eta^2(1-h)$ vanishes on $K$ quasi-everywhere, so the
Euler--Lagrange equation, first for Sobolev tests and then by
approximation, gives
\[
\begin{aligned}
&\int_{B(x,8R)}\eta^2|\nabla h|^2\,dy\\
&\quad=2\int_{B(x,8R)}\eta(1-h)\nabla h\cdot\nabla\eta\,dy + R^{-2}\int_{B(x,8R)}\eta^2h(1-h)\,dy\\
&\quad\le C_n MR^{n-2}.
\end{aligned}
\]
The Cauchy--Schwarz inequality and the annular energy estimate
bound the mixed term by $CMR^{n-2}$. The absorption term is at most {\color{black}$R^{-2}\int_{B(x,4R)}h\,dy\le CMR^{n-2}$}.
The function $\eta h$ is admissible for
$\operatorname{cap}_2(K,B(x,8R))$, and $M\le1$; hence
\begin{equation}\label{eq:appobstacleintermediate}
 \operatorname{cap}_2(K,B(x,8R))\le C_nR^{n-2}h(x).
\end{equation}
Extend $h$ by zero outside $B(x,8R)$. Comparison with $z_R$
in $G\cap B(x,8R)$ is valid because $(h-z_R)_+$ vanishes
quasi-everywhere on $E$ and on $\partial B(x,8R)$. Both satisfy \eqref{eq:appabsorption}
there. Testing their difference with $(h-z_R)_+$ gives
\[
 \int_{G\cap B(x,8R)}
 \bigl(|\nabla(h-z_R)_+|^2+R^{-2}(h-z_R)_+^2\bigr)\,dy=0.
\]
Thus $h\le z_R$. Taking increasing compact subsets of $E\cap B(x,R)$ and using
continuity from below of relative capacity proves
\eqref{eq:appobstacle}.

For $g\in C_c^\infty(G)$ and each $j\in\mathbb Z$, define
\[
 L_j=\{|g|\ge2^j\},\qquad
 h_j=\min\{1,2^{1-j}(|g|-2^{j-1})_+\}.
\]
Then $0\le h_j\le1$ and $h_j=1$ on $L_j$. Order
preservation gives $J_Rh_j\le J_R1$, hence
$h_j-J_Rh_j\ge z_R$ on $L_j$. Combining
\eqref{eq:appscale} with \eqref{eq:appobstacle} yields
\[
 \int_{L_j}\mathcal V_E\,dx
 \le C_n\int_G|\nabla h_j|^2\,dx.
\]
The disjoint strips $\{2^{j-1}<|g|<2^j\}$ give
\[
 \sum_j2^{2j}|\nabla h_j|^2\le4|\nabla g|^2,
 \qquad\sum_j2^{2j}\mathbf1_{L_j}\ge\tfrac13|g|^2.
\]
{\color{black}Multiplying the estimate for $\int_{L_j}\mathcal V_E\,dx$ by $2^{2j}$ and summing gives the bound for smooth $g$.} Density and Fatou's
lemma extend it to every $g\in H_0^1(G)$.
\end{proof}

\subsection{Sobolev traces and the forward implication}

Let $\Omega$ be a bounded $H^1$-extension domain, let
$F\subset\partial\Omega$ be closed, and fix a bounded extension
operator $E_\Omega:H^1(\Omega)\to H^1(\mathbb R^n)$.
By \cite[Theorem~1]{HKT08}, the interior measure density condition
holds with constants $b_0,r_0>0$ depending only on the extension
domain. For $y\in\partial\Omega$ and $0<r<r_0$, choose an interior
center within $r/2$ of $y$ and apply that condition at radius $r/2$.
Thus, after changing the constant, the boundary-centered estimate is
\begin{equation}\label{eq:appdensity}
 |\Omega\cap B(y,r)|\ge b r^n.
\end{equation}

For a function $U\in H^1(\mathbb R^n)$, its precise representative
$\widetilde U$ satisfies, outside a set of Sobolev $2$-capacity zero,
\begin{equation}\label{eq:appprecise}
 \lim_{r\downarrow0}\fint_{B(y,r)}
 |U(x)-\widetilde U(y)|\,dx=0.
\end{equation}
The inhomogeneous Sobolev capacity and the relative capacity in a
fixed containing ball have the same null sets. By
\eqref{eq:appdensity}, if $U_1$ and $U_2$ agree almost
everywhere in $\Omega$, their precise representatives agree
$2$-quasi-everywhere on $\partial\Omega$. Indeed, for
$W=U_1-U_2=0$ almost everywhere on $\Omega$,
\[
 b|\widetilde W(y)|\le r^{-n}\int_{B(y,r)}
 |W-\widetilde W(y)|\,dx\longrightarrow0
\]
by \eqref{eq:appprecise} and \eqref{eq:appdensity}.
The Havin--Bagby zero-trace criterion, in the formulation of
\cite[Proposition~1.1]{ET17}, therefore implies
\begin{equation}\label{eq:apptraceglobal}
 u\in H_F^1(\Omega)
 \quad\Longleftrightarrow\quad
 \widetilde{E_\Omega u}=0\text{ $2$-quasi-everywhere on }F.
\end{equation}
For the forward implication, extend the defining
approximants by $E_\Omega$. Their traces agree with those of the
original smooth functions by \eqref{eq:appdensity}, and convergence
in $H^1(\mathbb R^n)$ preserves zero trace. The reverse implication
follows by restricting the whole-space approximants from \cite{ET17}.

Choose a ball $B_*$ containing $\overline\Omega$ in its interior,
and $\chi\in C_c^\infty(B_*)$ equal to one near $\overline\Omega$.
If $u\in H_F^1(\Omega)$, \eqref{eq:apptraceglobal} gives
$\chi E_\Omega u\in H_0^1(B_*\setminus F)$.
The complement of $B_*\setminus F$ contains $F$; hence its
capacity weight dominates $\mathcal V_F$ on $\Omega$.
Lemma~\ref{lem:appfullhardy} implies
\begin{equation}\label{eq:appforward}
 \int_\Omega\mathcal V_F|u|^2\,dx
 \le C_n\int_{B_*\setminus F}|\nabla(\chi E_\Omega u)|^2\,dx
 \le C\|u\|_{H^1(\Omega)}^2.
\end{equation}
This proves \eqref{eq:A2capHardy} and one inclusion in
\eqref{eq:A2captrace}.

\subsection{A capacity-dominated measure and the converse}

\begin{lemma}[Measure dominated by relative capacity]
\label{lem:appmeasure}
If $K\Subset B_*$ is compact and
$\operatorname{cap}_2(K,B_*)>0$, there is a finite positive
Radon measure $\nu$ supported on $K$, with $\nu(K)>0$, such that, for every Borel set $A\subset B_*$,
\begin{equation}\label{eq:appmeasure}
 \nu(A)\le\operatorname{cap}_2(A,B_*).
\end{equation}
\end{lemma}
\begin{proof}
Let $h\in H_0^1(B_*)$ be the equilibrium potential of $K$.
Then $0\le h\le1$, $h=1$ quasi-everywhere on $K$, and
$\nu=-\Delta h$ is a positive Radon measure supported on $K$.
Its total mass is positive by the equilibrium identity.
For $\psi\in C_c^\infty(B_*)$, $\psi\ge0$, test the
equilibrium equation with $\psi^2/(h+\varepsilon)$.
The elementary inequality $2ab-a^2\le b^2$ gives
\[
 \int_{B_*}\frac{\psi^2}{h+\varepsilon}\,d\nu
 =\int_{B_*}\left(\frac{2\psi\nabla h\cdot\nabla\psi}{h+\varepsilon}
 -\frac{\psi^2|\nabla h|^2}{(h+\varepsilon)^2}\right)\,dx
 \le\int_{B_*}|\nabla\psi|^2\,dx.
\]
Since $h\le1$, taking $\varepsilon\downarrow0$, then the
infimum over admissible cutoffs for compact subsets, proves
\eqref{eq:appmeasure} on compact sets. Inner regularity of
$\nu$ and outer regularity of capacity give the Borel assertion.
This also shows $\nu$ charges no sets of capacity zero.
\end{proof}

\begin{lemma}[A one-scale estimate]\label{lem:appscaleconverse}
Let $u\in H^1(\Omega)$ and suppose a compact $K\subset F$ has
positive capacity and satisfies
$|\widetilde{E_\Omega u}|\ge2a$ on $K$, where $a>0$.
Let $\nu$ be as in Lemma~\ref{lem:appmeasure}. {\color{black}Put $r_k=2^{-k}r_0$ for integers $k\ge1$, where $r_0$ is the radius in the measure density condition preceding \eqref{eq:appdensity}. There are an integer $k_0\ge1$ with $r_{k_0}\le1$ and a measurable subset $K'\subset K$ with $\nu(K')>0$ such that, for every $k\ge k_0$,}
\begin{equation}\label{eq:appconverse1}
 \int_{K'}\kappa(y,r_k)\,d\nu(y)
 \le C a^{-2}\int_\Omega |u(x)|^2
     r_k^{-2} T_{r_k}(x)^2\,dx,
\end{equation}
where
\[
 \begin{aligned}
 \kappa(y,r)&=r^{2-n}\operatorname{cap}_2(F\cap B(y,r),B(y,8r)),\\
 T_r(x)&=r^{2-n}\operatorname{cap}_2(F\cap B(x,2r),B(x,8r)).
 \end{aligned}
\]
\end{lemma}
\begin{proof}
By \eqref{eq:appprecise}, \eqref{eq:appdensity}, and
$|\widetilde{E_\Omega u}|\ge2a$, for $\nu$-almost every
$y\in K$ and all sufficiently small $r$,
\begin{equation}\label{eq:apploweru}
 \int_{\Omega\cap B(y,r)}|u(x)|^2\,dx
 \ge c a^2 r^n.
\end{equation}
Indeed, the $L^1$ mean convergence and
\eqref{eq:appdensity} give
$\int_{\Omega\cap B(y,r)}|u|\,dx\ge
 a|\Omega\cap B(y,r)|\ge ab r^n$ for small $r$;
\eqref{eq:apploweru} follows by the Cauchy--Schwarz inequality.
{\color{black}Since \eqref{eq:apploweru} holds at all sufficiently small dyadic radii for $\nu$-almost every $y\in K$, choose an integer $k_0\ge1$ and a measurable subset $K'\subset K$ with $\nu(K')>0$ on which \eqref{eq:apploweru} holds for every $r=r_k$, $k\ge k_0$. Increase $k_0$, if necessary, so that $r_{k_0}\le1$ and $9r_{k_0}<\operatorname{dist}(K,\mathbb R^n\setminus B_*)$.}
By monotonicity and Lemma~\ref{lem:appcondenser}, for
$|x-y|<r$,
\[
 \kappa(y,r)\le C T_r(x).
\]
{\color{black}For $r=r_k$ with $k\ge k_0$, the choice of $k_0$ ensures that $B(x,8r)\subset B_*$ whenever $K'\cap B(x,r)$ is nonempty.} The measure
bound \eqref{eq:appmeasure} then gives
\[
 \begin{aligned}
 \nu(K'\cap B(x,r))
 &\le\operatorname{cap}_2(K'\cap B(x,r),B_*)\\
 &\le\operatorname{cap}_2(F\cap B(x,2r),B(x,8r))
 =r^{n-2}T_r(x).
 \end{aligned}
\]
The two bounds combine to give
$\int_{K'\cap B(x,r)}\kappa(y,r)\,d\nu(y)
 \le Cr^{n-2}T_r(x)^2$.
Integrate \eqref{eq:apploweru} against $\kappa(y,r)\,d\nu(y)$.
Tonelli's theorem then gives
\[
 \begin{split}
 \int_{K'}\kappa(y,r)\,d\nu(y)
 &\le Ca^{-2}r^{-n}\int_\Omega|u(x)|^2
     \int_{K'\cap B(x,r)}\kappa(y,r)\,d\nu(y)\,dx\\
 &\le Ca^{-2}r^{-2}\int_\Omega|u(x)|^2 T_r(x)^2\,dx.
 \end{split}
\]
This is \eqref{eq:appconverse1}.
\end{proof}

\begin{proof}[Completion of the proof of Proposition~\ref{prop:A2capacitary}]
Suppose $u\notin H_F^1(\Omega)$. By
\eqref{eq:apptraceglobal}, the precise representative of a
global extension is nonzero on a subset of $F$ of positive
capacity. Quasi-continuity and capacitability provide $a>0$
and a compact $K\subset F$ of positive relative capacity on
which $|\widetilde{E_\Omega u}|\ge2a$. Fix $\nu$ and $K'$
as in Lemmas~\ref{lem:appmeasure} and
\ref{lem:appscaleconverse}.

The Kellogg property \cite[Theorem~2]{HW83} gives divergence
of the Wiener integral for the Bessel capacity $\operatorname{Cap}_{1,2}$
at $2$-quasi-every $y\in F$. For $0<r\le1$ and
$A\subset B(y,r)$, zero extension and the Poincar\'e inequality give
$\operatorname{Cap}_{1,2}(A)\le C_n\operatorname{cap}_2(A,B(y,8r))$.
In particular, for $r_{k+1}=r_k/2$,
\[
 \int_{r_{k+1}}^{r_k}
 \frac{\operatorname{Cap}_{1,2}(F\cap B(y,t))}{t^{n-2}}\,
 \frac{dt}{t}\le C_n\kappa(y,r_k).
\]
Summing these inequalities gives
\begin{equation}\label{eq:appwiener}
 \sum_{k\ge k_0}\kappa(y,r_k)=\infty.
\end{equation}
This one-sided comparison remains valid for $n=2$.
Since $\nu$ charges no capacity-zero set,
\eqref{eq:appwiener} holds $\nu$-almost everywhere in $K'$.

For $2r\le R<4r$, monotonicity in the inner set and
Lemma~\ref{lem:appcondenser} yield
$\Theta_F(x,R)\ge c T_r(x)$, where $c$ depends only on $n$.
Since
$2\int_{2r}^{4r}R^{-3}\,dR=3/(16r^2)$, and the intervals
$[2r_k,4r_k)$ are disjoint,
\begin{equation}\label{eq:appscalesum}
 \sum_{k\ge k_0}r_k^{-2}T_{r_k}(x)^2
 \le C_n\mathcal V_F(x).
\end{equation}
Integrating \eqref{eq:appwiener} with respect to $\nu$
and then using \eqref{eq:appconverse1} and
\eqref{eq:appscalesum} gives
\[
 \infty
 =\sum_{k\ge k_0}\int_{K'}\kappa(y,r_k)\,d\nu(y)
 \le Ca^{-2}\int_\Omega\mathcal V_F(x)|u(x)|^2\,dx.
\]
{\color{black}Thus the finiteness of the weighted integral implies that $u\in H_F^1(\Omega)$.}
Together with \eqref{eq:appforward}, this proves both
assertions of Proposition~\ref{prop:A2capacitary}.
\end{proof}

\endgroup

\end{document}